\documentclass[twocolumn]{autart}    

\usepackage{cite}
\usepackage{amsfonts}
\usepackage{amsmath}
\usepackage{color}
\usepackage{textcomp}
\usepackage{graphics, graphicx}
\usepackage{enumerate}
\usepackage{array}
\usepackage{bm}
\usepackage{url} 
\usepackage{verbatim}
\usepackage{multirow}
\usepackage{mathrsfs}
\usepackage{algorithm}
\usepackage{amssymb}
\usepackage{stfloats}
\usepackage{algpseudocode}
\usepackage{booktabs}

\usepackage{apacite} 

\newtheorem{mytheo}{Theorem}
\newtheorem{mydef}{Definition}

\newtheorem{mycoro}{Corollary}

\newtheorem{mylemma}{Lemma}

\newtheorem{myassump}{Assumption}
\newtheorem{myprop}{Proposition}
\newtheorem{myexa}{Example}
\newtheorem{myrmk}{Remark}

\usepackage{bm}

\edef\endfrontmatter{%
  \unexpanded\expandafter{\endfrontmatter}
  \noexpand\endNoHyper 
}  

\begin{document}

\begin{frontmatter}

\title{Chance-constrained probability measure optimization\thanksref{footnoteinfo}
} 

\thanks[footnoteinfo]{This paper was not presented at any IFAC 
meeting. Corresponding author Y. Wang. }

\author[TUAT]{Xun Shen}\ead{xun.shen.whu@gmail.com},    
\author[Melbourne]{Ye Wang}\ead{ye.wang@unimelb.edu.au},
\author[DUT]{Yuhu Wu}\ead{wuyuhu@dlut.edu.cn},               
\author[ISM]{Satoshi Ito}\ead{sito@ism.ac.jp}, 
\author[IST]{Jun-ichi Imura}\ead{imura@sc.e.titech.ac.jp}  

\address[TUAT]{Graduate School of Engineering, Tokyo University of Agriculture and Technology, Tokyo, 184-8588 Japan}  
\address[Melbourne]{School of Mathematics and Statistics, The University of Melbourne, VIC, 3010, Australia}
\address[DUT]{School of Control Science and Engineering, Dalian University of Technology, Dalian, 116024, China}        
\address[ISM]{Department of Statistical Inference and Mathematics, The Institute of Statistical Mathematics, Tokyo, 190-8562, Japan}             
\address[IST]{Department of Systems and Control Engineering, Institute of Science Tokyo, Tokyo, 152-8550, Japan}        


\begin{abstract}                          
Stochastic optimization with chance constraints often relies on deterministic decision-making, where optimal decisions are fixed and may be applied to a system repeatedly. A critical question arises: Can probabilistic decision-making, where probability measures/distributions are considered as decision variables, outperform deterministic decision-making in terms of the expected performance when chance constraints are present? This paper addresses this question by introducing the \textit{Chance-Constrained Probability Measure Optimization} (CCPMO) framework, which formulates the problem of optimizing probabilistic decisions under chance constraints. We first establish the existence of the optimal solution to CCPMO. Crucially, we prove that the optimal probabilistic decisions can always be represented by a probability measure concentrated on only two points, thereby reducing the CCPMO problem to an equivalent, simplified form. To solve this reduced problem, we propose a sample-based smooth approximation method. This approach leverages samples of model uncertainties to construct an approximate problem with uniform convergence and probabilistic feasibility guarantees. The approximate problem can be efficiently solved using standard nonlinear programming techniques. Finally, we validate the proposed framework through a numerical example of a quadrotor control problem under turbulent conditions. The results demonstrate that probabilistic decision-making outperforms deterministic approaches in expected performance while satisfying safety-critical chance constraints.
\end{abstract}

\end{frontmatter}

\section{Introduction}\label{intro}

Robust decision-making against uncertainty is often achieved by solving an optimization problem with chance constraints that ensure solutions satisfying probabilistic safety requirements \cite{Liu:2023}. Recent advances in this domain span Reinforcement Learning (RL) and Stochastic Model Predictive Control (SMPC).  For example, constrained RL frameworks enforce expected risk-return constraints during policy updates, enabling safe exploration in safety-critical applications \cite{Wachi, Yan2025OfflineGuarded}. These methods train neural networks to produce optimal decisions that adhere to predefined constraints, provided sufficient training data. Conversely, SMPC recursively solves finite-horizon optimal control problems with chance constraints to obtain an implicit safe policy \cite{Mayne:2016}. However, SMPC relies on precise system models, limiting its applicability in data-driven settings. Bridging this gap, recent work integrates SMPC with RL, where RL learns parameters for the objective function, model and constraints of SMPC, while SMPC guides safe exploration during training \cite{Gros}.

Existing optimization methods predominantly focus on \emph{deterministic decision-making}, where decisions are selected deterministically. This approach can be treated as a special case of \emph{probabilistic decision-making}, where decisions are selected probabilistically.  While probabilistic decisions have been used in data-driven SMPC and RL for exploration \cite{Bravo}, their potential to systematically outperform deterministic decisions under chance constraints remains unexplored. Furthermore, the problem of formally optimizing probabilistic decisions under chance constraints has yet to be rigorously addressed. This paper formalizes \emph{Chance-Constrained Probability Measure Optimization} (CCPMO), which optimizes probabilistic decisions (probability measure over decision variables) under chance constraints. 

Chance-constrained optimization for deterministic decision-making with finite-dimensional decision variables has been extensively studied in stochastic programming and control engineering. However, it is generally NP-hard due to non-convex feasible sets \cite{Shapiro}. Current research mainly follows two streams: (a) impose specific structure for constraint functions or the distribution of random variables, such as linear or convex constraint functions \cite{Nemirovski}, finite sample space of random variables \cite{Luedtke:2008}, elliptically symmetric Gaussian-similar distributions \cite{Ackooij}; (b) extract samples to approximate the chance constraints, i.e. sample-based methods \cite{Campi:2008, Pena:2020, Geletu:2019}. The most common sample-based method is the scenario approach \cite{Campi:2008, Campi_unconvex}, which approximates the original problem by generating a deterministic problem from samples, with feasibility improving as sample size increases. The sample-average approach \cite{Luedtke:2008, Geletu:2019, Pena:2020} ensures both feasibility and optimality. However, these approaches are developed for deterministic decision variables with a finite-dimensional parameterization, whereas CCPMO operates in an infinite-dimensional probability measure space and optimizes over probability measures to produce probabilistic decisions.

Optimization with an infinite number of chance constraints in finite-dimensional decision variable space has also been intensively studied,\cite{Ackooij:2020, Berthold:2022}. In \cite{Ackooij:2020}, the variational tools are applied to formulate generalized differentiation of chance/robust constraints. The method of getting the explicit outer estimations of sub-differentials from data is also established. An adaptive grid refinement algorithm is developed to solve the optimization with chance/robust constraints in \cite{Berthold:2022}. In addition, probability functions and their (generalized) derivatives have been studied for possibly infinite-dimensional deterministic decision variables in Banach spaces under Gaussian distributions, e.g., \cite{Hantoute2019}.

Chance constraints combined with recourse and multi-stage decision making have also been widely studied, where the decision remains a deterministic first-stage action (and possibly deterministic recourse decisions or decision rules), rather than a probability measure \cite{Liu2016}. Likewise, dynamic chance-constraint formulations in infinite-dimensional function spaces typically optimize deterministic policies/decision rules (e.g., continuous functions over compact sets or time intervals), rather than optimizing a probability measure as the decision variable \cite{GonzalezGradon2020, Geletu:2020, Grandon:2022}.

The above research on optimization with chance constraints in finite-dimensional vector spaces can still only prove convergence in the case where the dimension of the decision variable is infinite. In \cite{Cui:2022}, a new formulation of a stochastic program called Affine Chance Constrained Stochastic Programs (ACC-SP) is presented. The constraint of ACC-SP includes mixed-signed affine combinations of countable violation probabilities. Our addressed CCPMO problem differs from ACC-SP in the following two points: (1) CCPMO optimizes a probability measure defined on the $\sigma$-algebra of Borel sets, which ACC-SP optimizes a decision variable on a finite-dimensional set; (2) The constraint of CCPMO is an integral of violation probabilities with respect to the probability measure to be optimized. While the constraint of ACC-SP involves fixed coefficients of affine combinations. Besides, CCPMO aims to optimize probabilistic decisions, while ACC-SP focuses on deterministic decisions.

The main contribution of this paper is to propose a novel CCPMO framework with theoretical guarantees for probabilistic decision-making. The specific contributions are summarized as follows: (a) Formalization of CCPMO on Borel probability measure space, establishing a foundation for probabilistic decision-making under chance constraints (Section \ref{sec:problem}); (b) Theoretical guarantees on the existence of optimal CCPMO solution under a mild assumption (Section~\ref{sec:existence}); (c) CCPMO reduction with finite-dimensional decision variables, based on the optimal solution of CCPMO that has a probability measure concentrated on two decisions (Section~\ref{sec:reduction}); (d) Sample-based smooth approximation with uniform convergence and probabilistic feasibility guarantees (Section~\ref{sec:approximation}).

\section{Problem Formulation}\label{sec:problem}

Consider the following chance-constrained optimization: Let $x\in\mathcal{X}\subset\mathbb{R}^n$ be the decision variable, where $\mathcal{X}$ is a compact set. The objective function $J:\mathcal{X}\rightarrow\mathbb{R}$ is a continuous scalar function on $\mathcal{X}$. The optimization problem is further subject to a chance constraint, which ensures that the probability of satisfying a given condition exceeds a specified threshold $1-\alpha$ with $\alpha \in [0,1]$. 
Let $\Omega$ be a sample space and $\mathcal{F}$ be the $\sigma$-algebra of Borel sets of $\Omega$. Let $\omega\in\Omega$ be a sample in $\Omega$. With a probability measure $\mathsf{Pr}_{\omega}$ defined on Borel space $\left(\Omega,\mathcal{F}\right)$, a probability space $\left(\Omega,\mathcal{F},\mathsf{Pr}_{\omega}\right)$ is formed. The constraint function involves a random variable $\xi:\Omega\rightarrow\Xi\subseteq\mathbb{R}^s$. Let $\Xi_{\mathsf{Leb}}\subseteq\Xi$ be a subset of $\Xi$. We define the probability of having $\xi\in\Xi_{\mathsf{Leb}}$ by $\mathsf{Pr}_{\xi}\left\{\xi\in \Xi_{\mathsf{Leb}}\right\}:=\displaystyle\int_{\xi^{-1}\left(\Xi_{\mathsf{Leb}}\right)} \mathsf{d}\mathsf{Pr}_{\omega}\left(\omega\right)$. Assume that there is a continuous probability density function $p(\xi)$ with the support $\Xi$. Then, we can rewrite $\mathsf{Pr}_{\xi}\left\{\xi\in \Xi_{\mathsf{Leb}}\right\}$ by $\mathsf{Pr}_{\xi}\{\xi\in\Xi_{\mathsf{Leb}}\}:=\displaystyle\int_{\Xi_{\mathsf{Leb}}}p(\xi)\mathsf{d}\xi$.

The constraint function $h:\mathcal{X}\times\Xi\rightarrow\mathbb{R}^m$ is a continuous vector-valued function with $h(x,\xi):=[h_1(x,\xi),...,h_m(x,\xi)]^\top$. 
Define $\bar{h}:\mathcal{X}\times\Xi\rightarrow\mathbb{R}$ by 
\begin{equation}\label{eq:def_bar_h}
    \bar{h}(x,\xi) :=\max_{i=1,\ldots,m} h_i(x,\xi). 
\end{equation}
Let $\mathbb{I}\{y\}$ present the indicator function with $\mathbb{I}\{y\}=1$ if $y\leq 0$ and $\mathbb{I}\{y\}=0$ otherwise. 
Furthermore, we define $\mathbb{P}:\mathcal{X}\rightarrow[0,1]$ by 
\begin{equation}
\label{eq:def_P_x}
    \mathbb{P}(x):=\int_{\Xi}\mathbb{I}\{\bar{h}(x,\xi)\}p(\xi)\mathsf{d}\xi.
\end{equation}
Note that $\mathbb{P}(x)$ is the probability that $\bar{h}(x,\xi)\leq 0$ holds for a given $x$. 
Therefore, the chance constraint is given by
\begin{equation}\label{eq:chance constraint}
    \mathbb{P}(x)\geq 1-\alpha,
\end{equation}
where $\alpha \in [0, 1]$ is a violation probability threshold.

In this work, after the chance-constrained optimization is solved, the resulting optimal decision $x^*$ will be implemented multiple times. To provide context, we first review deterministic decisions for the considered chance-constrained optimization. Since the decision $x^*$ will be applied repeatedly, we then shift our focus to probabilistic decision-making. Specifically, we address a chance-constrained probability measure optimization problem, which generalizes the deterministic case with improved performance and accounts for the repeated implementation of decisions.

\subsection{Deterministic decision-making}

According to the above setup, the chance-constrained optimization problem for deterministic decision-making can be formulated as
\begin{align} 
\label{eq:Q_alpha}
&\min_{x\in \mathcal{X}} \,\, J(x) \tag{$Q_\alpha$}, \\
&{\normalfont \mathsf{s.t.}}\quad  \mathbb{P}(x)\geq 1-\alpha. \nonumber
\end{align} 

Let the feasible set of Problem \ref{eq:Q_alpha} be
\begin{equation}\label{eq:mathX_alpha_Q}
    \mathcal{X}_{\alpha}:=\{x\in\mathcal{X}:\mathbb{P}(x)\geq 1-\alpha\},
\end{equation}
and the optimal objective value $J^*_\alpha$ and optimal
solution set $X_\alpha$ of \ref{eq:Q_alpha} be
\begin{align}
    \label{eq:J_star_alpha_Q}
    J^*_{\alpha} &:= \min\{J(x):x\in\mathcal{X}_{\alpha}\},\\
    \label{eq:X_opt_alpha_Q}
    X_{\alpha} &:= \{x\in\mathcal{X}_{\alpha}:J(x)=J^*_{\alpha}\},    
\end{align}
respectively. 

Observe that if there exists $\underline{\alpha}\in[0,1]$ such that $\mathcal{X}_{\underline{\alpha}}\neq\emptyset$,
then $\mathcal{X}_{\alpha}\neq\emptyset$ for all $\alpha\in[\underline{\alpha},1]$.
Hence, throughout the paper we restrict attention to violation levels $\alpha\in[\underline{\alpha},1]$,
and assume that for every such $\alpha$ the problem $Q_\alpha$ is feasible and attains a global minimum with finite optimal value.

By solving Problem \ref{eq:Q_alpha}, we can obtain an optimal deterministic decision $x^*\in X_\alpha$. As aforementioned, the deterministic decision $x^*$ will be implemented multiple times. Through the following example, we show a potential improvement on performance using probabilistic decisions compared to deterministic one.

\begin{myexa}\label{exam1} 
The set $\mathcal{X}$ is given by $\mathcal{X}=[-2,2]\subset \mathbb{R}$.
The objective function $J(x)$ and the single constraint function $h(x,\xi)$ are given as follows:
\begin{subequations}
    \begin{align}
       J(x) &= -(x+0.6)^2+2, \label{eq:J_num_1} \\
       h(x,\xi) &= x-1.4+\xi, \label{eq:h_num_1}
    \end{align}
\end{subequations}
where $\xi\sim\mathcal{N}(m_{\xi},\Sigma_{\xi})$ with $m_{\xi}=0$, and $\Sigma_{\xi}=1$. The violation probability threshold $\alpha$ is set as $\alpha = 0.25$. The profiles of the objective function $J(x)$ and the violation probability $1-\mathbb{P}(x)$ are plotted in Fig. \ref{fig:exa_one} (a) and (b), respectively. It can be seen that as $x$ increases, the violation probability $1-\mathbb{P}(x)$ increases while the objective function $J(x)$ decreases. Therefore, the optimal solution $x^*_{\alpha}$ is at the boundary with $1-\mathbb{P}(x^*_{\alpha})=0.25$. Since $\xi$ is a random variable that obeys normal distribution and $h(x,\xi)$ is linear, we can obtain $x^*_{\alpha}=0.7259$ from checking the standard normal distribution table, and the
corresponding optimal objective value is $J^*_{\alpha}=0.2420$.

\begin{figure}[h]
\centering
\includegraphics[scale=0.375]{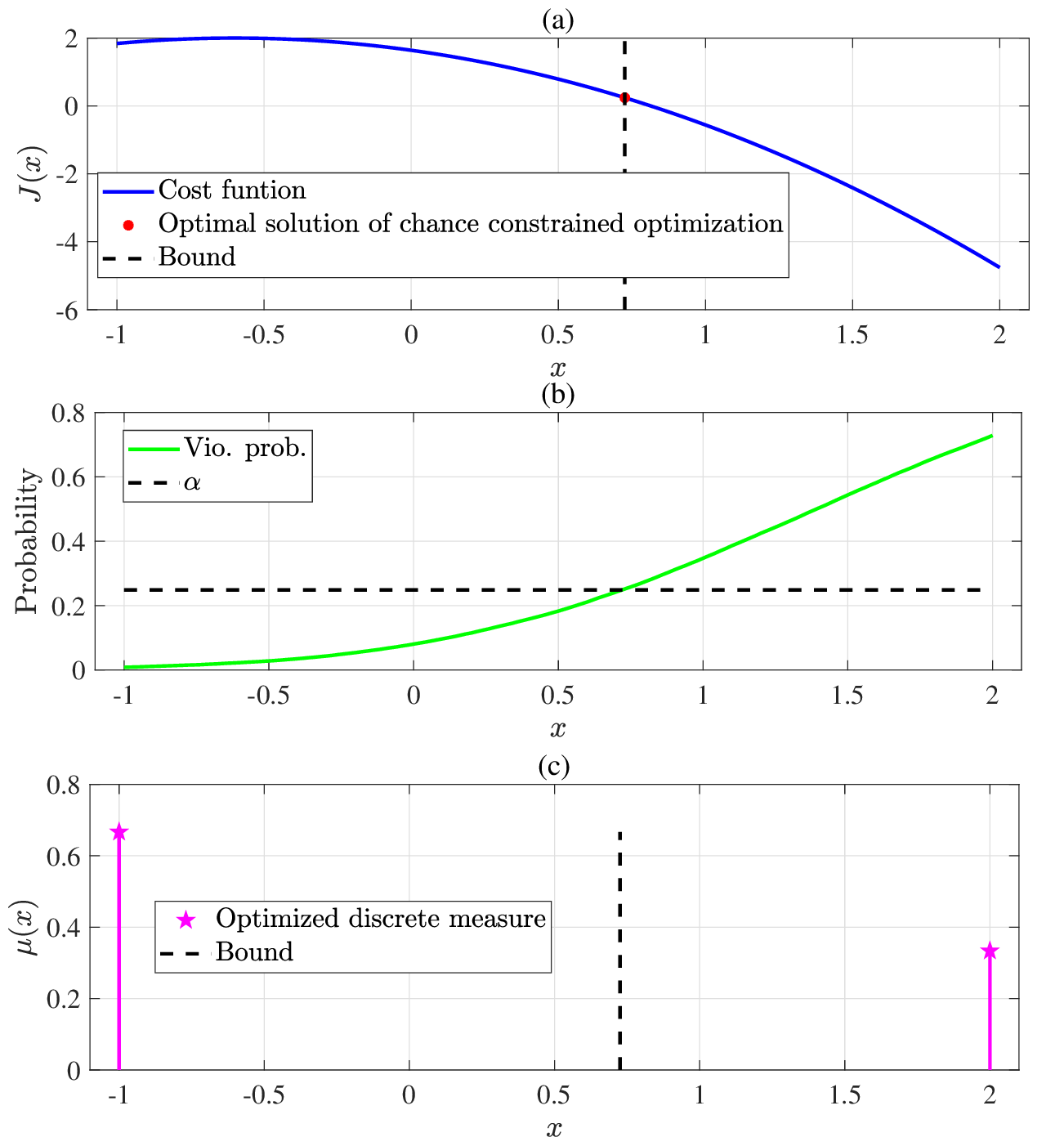}
\centering
\caption{Example in one dimension: (a) Profile of $J(x)$ and optimal solution of chance-constrained optimization; (b) Profile of $1-\mathbb{P}(x)$; (c) Optimized discrete measure.}
\label{fig:exa_one}
\end{figure}

We then show how probabilistic decisions can further improve the expected objective value. Let $\mathcal{C}_2=\{x^{(1)},x^{(2)}\}$ be a set with two samples from
$\mathcal{X}=[-2, 2]$, where $x^{(1)}=-1$ and $x^{(2)}=2$. 
The objective function values of these two samples are $J(x^{(1)})=1.84$ and $J(x^{(2)})=-4.76$. 
The violation probabilities are $\mathbb{P}_{\mathsf{vio}}(x^{(1)})=1-\mathbb{P}(x^{(1)})=0.008$ and $\mathbb{P}_{\mathsf{vio}}(x^{(2)})=1-\mathbb{P}(x^{(2)})=0.729$, respectively. Given a discrete probability measure $\mu$ on $\mathcal{X}=[-2, 2]$ with $\mu\left(\{x^{(1)}\}\right)=0.67$ and $ \mu\left(\{x^{(2)}\}\right)=0.33$.
Using this measure, $x^{(1)}$ is chosen with probability $0.67$ and $x^{(2)}$ is chosen with probability $0.33$. Then, the expected objective value is 
\begin{equation*}
    J(x^{(1)})\mu(\{x^{(1)}\})+J(x^{(2)})\mu(\{x^{(2)}\})=-0.338<J^*_\alpha,
\end{equation*}
and the expected violation probability is
\begin{equation*}
    \mathbb{P}_{\mathsf{vio}}(x^{(1)})\mu(\{x^{(1)}\})+\mathbb{P}_{\mathsf{vio}}(x^{(2)})\mu(\{x^{(2)}\})=0.2459<\alpha,
\end{equation*}
which implies a  probabilistic decision that satisfies the chance constraint can improve the expected control performance compared with the deterministic optimal decision
 $x^*_{\alpha}$.
\end{myexa}

\begin{myrmk}
    Example \ref{exam1} demonstrates that optimized discrete measure on two random $x^{(1)}$ and $x^{(2)}$ is a feasible probabilistic decision and can outperform the optimal deterministic decision in terms of optimal objective value meanwhile maintaining chance constraint.
\end{myrmk}

\subsection{Probabilistic decision-making and CCPMO}\label{subsec:prob_form}

We now formally formulate the optimization problem to identify the optimal probabilistic decision. Unlike deterministic optimization, which directly optimizes the decision variable $x$, the probabilistic problem seeks to optimize random variable $x$ characterized by a probability measure $\mu$.

Let $\mathscr{B}(\mathcal{X})$ be the Borel $\sigma$-algebra ($\sigma$-field) on $\mathcal{X} \subset \mathbb{R}^n $ with Euclidean distance. Notice that $\left(\mathcal{X},\mathscr{B}(\mathcal{X})\right)$ is a measurable space (Borel space). 
 In this work, rather than 
 fixing one measure  $\mu$ on $\left(\mathcal{X},\mathscr{B}(\mathcal{X})\right)$, it is interested to consider
the totality of possible (probability) measures. 
Hence, we denote by   $M(\mathcal{X})$ as  the set of all Borel probability measures on  $\left(\mathcal{X},\mathscr{B}(\mathcal{X})\right)$. 

Denote the space of continuous $\mathbb{R}$-valued functions on the compact set $\mathcal{X}$ by
\begin{equation}
\label{eq:C_X_R}
     \mathscr{C}(\mathcal{X}):=\{f:\mathcal{X}\rightarrow\mathbb{R} \, |f\ \text{is continuous}\}.
\end{equation}
It is able to define a metric on $\mathscr{C}(\mathcal{X})$ by  
\begin{equation}
\label{eq:tau}
    \tau(f,f'):=\|f-f'\|_\infty,\ \forall f,f'\in\mathscr{C}(x),
\end{equation}
where $\|f\|_\infty$ is defined as $\|f\|_\infty:=\sup_{x\in\mathcal{X}}|f(x)|.$

The metric $\tau(\cdot,\cdot)$ turns $\mathscr{C}(\mathcal{X})$ into a complete metric space. 
Define the pseudo-distance (pseudo-metric)  between any two probability measures $\mu,\nu\in M(\mathcal{X})$ associated with $f\in\mathscr{C}(\mathcal{X})$ by 
\begin{equation}
    \label{eq:distance_mu_nu}
    \mathcal{W}_f(\mu,\nu)=\left|\int_{\mathcal{X}}f(x)\mathsf{d}\mu-\int_{\mathcal{X}}f(x)\mathsf{d}\nu\right|.
\end{equation}
The weak$^*$ convergence of probability measures is defined as follows \cite{Billingsley}.

\begin{mydef}
    \label{def:weak_convergence}
    Let $\{\mu_k\}_{k=0}^{\infty}$ be a sequence in $M(\mathcal{X})$. We say that $\{\mu_k\}_{k=0}^{\infty}$ converges weakly$^*$ to $\mu$ if
    \begin{equation}
        \label{eq:weak_converge_def}        \lim_{k\rightarrow\infty}\mathcal{W}_f(\mu_k,\mu)=0,\ \forall f\in\mathscr{C}(\mathcal{X}).
    \end{equation}
    Besides, we say that $\{\mu_k\}_{k=0}^{\infty}$ converges to $\mu$ associated with $f\in\mathscr{C}(\mathcal{X})$ if 
    \begin{equation}
        \label{eq:converge_f_def}        \lim_{k\rightarrow\infty}\mathcal{W}_f(\mu_k,\mu)=0.
    \end{equation}
\end{mydef}
Notice that $h(x,\xi)$ is a random variable for a given $x\in\mathcal{X}$. We introduce an assumption on $h(x,\xi)$ for the continuity analysis of $\mathbb{P}(x)$. Let $\mathsf{cl}\{Z\}$ be the closure of a set $Z$. Let $\mathsf{supp}\ p:=\mathsf{cl}\{\xi\in\Xi:p(\xi)>0\}$, and, for each $x\in\mathcal{X}$, $\Xi^{\mathsf{supp}}(x):=\{\xi\in\mathsf{supp}\ p:\bar{h}(x,\xi)=0\}$. The following assumption holds throughout the paper.

\begin{myassump}
\label{assump:J_h}
For each $x\in\mathcal{X}$, $\mathsf{Pr}_\xi\left\{\Xi^{\mathsf{supp}}(x)\right\}=0$ holds.
\end{myassump}

Assumption \ref{assump:J_h} implies the continuity of $\mathbb{P}(x)$ \cite{Geletu:2019, Kibzun}, which guarantees  that the feasible set $\mathcal{X}_\alpha$ defined by \eqref{eq:mathX_alpha_Q} and optimal solution set $X_\alpha$ given by \eqref{eq:X_opt_alpha_Q} are measurable sets for all $\alpha\in[0,1]$.

With the above notation and discussion, a CCPMO Problem \ref{eq:P_alpha} associated with Problem \ref{eq:Q_alpha} can be formulated as follows:
\begin{align} 
    \label{eq:P_alpha}
    &\underset{\mu\in M(\mathcal{X})}{{\normalfont\mathsf{min}}} \,\, \int_{\mathcal{X}}J(x)\mathsf{d}\mu, \tag{$P_\alpha$} \\
    &{\normalfont \mathsf{s.t.}}\quad  \int_{\mathcal{X}}\mathbb{P}(x) \mathsf{d}\mu \geq 1-\alpha, \nonumber
\end{align} 
where $\mathbb{P}(x)$ is defined in \eqref{eq:def_P_x}. 

At this stage, the notation ‘$\mathsf{min}$’ is formal. 
The attained optimum is established later in Theorem \ref{theo:P_alpha_opt_solution}. Equivalently, one may read it as ‘$\mathsf{inf}$’ before Theorem \ref{theo:P_alpha_opt_solution}.

In Problem \ref{eq:P_alpha}, the uncertainty of $h(x,\xi)$ comes from both $x$ and $\xi$, since $x$ obeys the distribution dominated by probability measure $\mu$. Let $M_{\alpha}(\mathcal{X})$ be the feasible set of Problem \ref{eq:P_alpha} by
\begin{equation}
\label{eq:M_alpha}
    M_{\alpha}(\mathcal{X}):=\Big\{\mu\in M(\mathcal{X}):\int_{\mathcal{X}}\mathbb{P}(x) \mathsf{d}\mu \geq 1-\alpha\Big\}.
\end{equation}
The optimal objective value and solution set of Problem \ref{eq:P_alpha} are written as 
\begin{align}
\label{eq:mathJ_star}
    \mathcal{J}^*_{\alpha}&:=\inf_{\mu\in M_{\alpha}(\mathcal{X})} \int_{\mathcal{X}}J(x) \mathsf{d}\mu, \\
\label{eq:A_alpha}
    A_{\alpha}&:=\Big\{\mu\in M_{\alpha}(\mathcal{X}):\int_{\mathcal{X}}J(x) \mathsf{d}\mu = \mathcal{J}^*_{\alpha}\Big\}.    
\end{align}
A probability measure $\mu^*_{\alpha}\in A_{\alpha}$ is called an optimal probability measure for Problem \ref{eq:P_alpha}. The optimal probabilistic decision generates a decision at every trial from $\mathcal{X}$ probabilistically according to the optimal probability measure $\mu^*_\alpha$. If the decision process repeats many times, the mean value of the objective value is minimal.

\section{Existence of the optimal CCPMO solution}
\label{sec:existence}

In this section, we show the existence of at least an optimal solution to the CCPMO Problem \ref{eq:P_alpha}. Additionally, we present a proposition that establishes a relationship between the optimal objective values of Problems \ref{eq:P_alpha} and \ref{eq:Q_alpha}.

\begin{mytheo}\label{theo:P_alpha_opt_solution}
If Assumption \ref{assump:J_h} holds and $\mathcal{X}$ is compact, we have $A_\alpha\neq\emptyset$.
\end{mytheo}

\begin{pf}
    The optimal objective value $\mathcal{J}^*_\alpha$ of Problem \ref{eq:P_alpha} is defined by \eqref{eq:mathJ_star}. We aim to prove that there exists at least one $\mu^*\in M_\alpha(\mathcal{X})$ such that $\displaystyle\int_{\mathcal{X}}J(x) \mathsf{d}\mu^*= \mathcal{J}^*_{\alpha}$. 
    
    For any $\mu \in M(\mathcal{X})$, we have
    $$
    \int_{\mathcal{X}}J(x)\mathsf{d}\mu\ge 
    \int_{\mathcal{X}} \min_{x\in \mathcal{X}} \,\, J(x) \mathsf{d}\mu=\min_{x\in \mathcal{X}} \,\, J(x),
    $$
    which implies that $\mathcal{J}^*_{\alpha}>-\infty$. Then, there exists a sequence $\{\mu_{k}\}_{k=1}^{\infty}\subset M_\alpha(\mathcal{X})$ such that $\displaystyle \int_{\mathcal{X}}J(x)\mathsf{d}\mu_k\rightarrow\mathcal{J}^*_{\alpha}$. 
    Since $\mathcal{X}$ is assumed to be compact, by Prokhorov Theorem (Theorem 8.6.2 of \cite{Bogachev}), there exists a subsequence $\{\mu_{k_s}\}_{s=1}^{\infty}\subset\{\mu_k\}_{k=1}^\infty$ such that converges weakly$^*$ to $\mu^*\in M(\mathcal{X})$ in the sense of Definition \ref{def:weak_convergence}. Since $\displaystyle \int_{\mathcal{X}}J(x)\mathsf{d}\mu_{s_k} \rightarrow \displaystyle \int_{\mathcal{X}}J(x)\mathsf{d}\mu^*$ by $J(\cdot)\in \mathscr{C}(\mathcal{X})$, we have
    \begin{equation}
        \label{eq:mu_star_optimality_proof}
        \int_{\mathcal{X}} J(x)\mathsf{d}\mu^*
        =\lim_{s\to \infty}\int_{\mathcal{X}} J(x)\mathsf{d}\mu_{n_s}=
        \mathcal{J}^*_{\alpha}.
    \end{equation}
    
    By Assumption \ref{assump:J_h}, we know $\mathbb{P}(\cdot)\in\mathscr{C}(\mathcal{X})$.
    Together with $\{\mu_{k}\}_{k=1}^{\infty}\subset M_\alpha(\mathcal{X})$, we have that  $\displaystyle \int_{\mathcal{X}}\mathbb{P}(x)\mathsf{d}\mu_{k_s}\geq 1-\alpha$, $ \forall s\ge 1$. As a result, since $\mu_{k_s}$ converges weakly$^*$ to $\mu^*$, 
    \begin{equation}
        \label{eq:mu_star_feasibility_proof}
        \int_{\mathcal{X}}\mathbb{P}(x)\mathsf{d}\mu^*(x)= \lim_{s\to \infty}\int_{\mathcal{X}}\mathbb{P}(x)\mathsf{d}\mu_{k_s}(x)\geq 1-\alpha,
    \end{equation}
    which implies $\mu^*\in M_\alpha(\mathcal{X})$.
    By \eqref{eq:mu_star_optimality_proof} and \eqref{eq:mu_star_feasibility_proof}, $\mu^*$ is an optimal solution of Problem \ref{eq:P_alpha}. Thus, $A_\alpha\neq\emptyset$ holds. \qed
\end{pf}

Theorem \ref{theo:P_alpha_opt_solution} establishes that the optimal solution set $A_\alpha$ of Problem \ref{eq:P_alpha} is non-empty. Therefore, the infimum ``$\inf$" in \eqref{eq:mathJ_star} can be replaced with the minimum ``$\min$", as the existence of an optimal solution guarantees the attainment of the infimum.

We now present the following proposition, which establishes the relationship between $\mathcal{J}^*_\alpha$ and $J^*_\alpha$.

\begin{myprop}
    \label{prop:calJ_J_star}
    Suppose that Assumption \ref{assump:J_h} holds. Then, $$\mathcal{J}^*_\alpha\leq J^*_\alpha,$$ where $\mathcal{J}^*_\alpha$ and $J^*_\alpha$ are the optimal objective values of Problems \ref{eq:P_alpha} and \ref{eq:Q_alpha}, respectively.
\end{myprop}
\begin{pf}
    By definition, the chance-constrained feasible set of Problem \ref{eq:Q_alpha} is
    $\mathcal{X}_\alpha:=\{x\in\mathcal{X}:\mathbb{P}(x)\geq 1-\alpha\}$.
    Under Assumption \ref{assump:J_h}, $\mathbb{P}(\cdot)$ is continuous on $\mathcal{X}$ .
    Hence, $\mathcal{X}_\alpha$ is the inverse image of the closed set $[1-\alpha,1]$ under a continuous mapping, which implies that $\mathcal{X}_\alpha$ is closed in $\mathcal{X}$.
    Therefore, $\mathcal{X}_\alpha$ is a Borel set, i.e., $\mathcal{X}_\alpha\in\mathcal{F}$.
    Moreover, since $J(\cdot)$ is continuous on $\mathcal{X}$ by Assumption \ref{assump:J_h}, the optimal solution set $X_\alpha$, defined by \eqref{eq:X_opt_alpha_Q}, is a closed subset of $\mathcal{X}_\alpha$ and thus also Borel measurable.
    
    Since $X_\alpha\neq\emptyset$ for the violation levels considered in this paper, pick any $\tilde{x}_\alpha\in X_\alpha$ and define the Dirac probability measure $\tilde{\mu}_\alpha:=\delta_{\tilde{x}_\alpha}\in M(\mathcal{X})$. Then $\tilde{\mu}_{\alpha}(X_{\alpha})=1$.
    Then, we have
    \begin{align}  
    \label{eq:J_m_mu}
    \int_{\mathcal{X}}J(x)\mathsf{d}\tilde{\mu}_{\alpha}&=\int_{X_{\alpha}}J(x)\mathsf{d}\tilde{\mu}_{\alpha}=J^*_{\alpha},\\
    \label{eq:P_m_mu}
    \int_{\mathcal{X}}\mathbb{P}(x)\mathsf{d}\tilde{\mu}_{\alpha}&=\int_{X_{\alpha}}\mathbb{P}(x) \mathsf{d}\tilde{\mu}_{\alpha}\geq 1-\alpha.   
    \end{align} 
    Thus, $\tilde{\mu}_{\alpha}$ is a feasible solution for \ref{eq:P_alpha} with objective value as $J^*_{\alpha}$, namely, $\tilde{\mu}_{\alpha}\in M_{\alpha}(\mathcal{X})$ and $\mathcal{J}^*_\alpha\leq J^*_\alpha$. \qed
\end{pf}

Proposition \ref{prop:calJ_J_star} demonstrates that the optimal objective value $\mathcal{J}^*_\alpha$ of Problem \ref{eq:P_alpha} may be smaller than the optimal objective value $J^*_\alpha$ of Problem \ref{eq:Q_alpha}. This raises two underlying critical questions: (1) Under what conditions does $\mathcal{J}^*_{\alpha} = J^*_{\alpha}$? (2) How can Problem \ref{eq:P_alpha} be solved when $ \mathcal{J}^*_{\alpha} < J^*_{\alpha} $?
To address the first question, we present Corollary \ref{coro:W_alpha_Q_alpha} in Section \ref{sec:reduction}, which provides two sufficient conditions for $\mathcal{J}^*_{\alpha} = J^*_{\alpha}$. For the second question, we first derive an equivalent reduction of Problem \ref{eq:P_alpha} using Theorem~\ref{theo:P_alpha_2_point_approx} proposed in Section \ref{sec:reduction}. Subsequently, in Section~\ref{sec:approximation}, we propose a sample-based smooth approximation to solve the reduced problem.

\section{CCPMO Problem Reduction}\label{sec:reduction}

In this section, we show the existence of an optimal probabilistic decision for Problem \ref{eq:P_alpha} whose probability measure is concentrated on two points. This result leads to an equivalently reduced form of Problem \ref{eq:P_alpha}.

We give an assumption on the optimal solution set $A_\alpha$. 

\begin{myassump}\label{assump:P_alpha_opt_sol_seq}
    There exists an optimal solution $\mu^*\in A_\alpha$ of Problem \ref{eq:P_alpha} such that for any $\delta>0$ there exists a probability measure $\mu$, different from $\mu^*$, such that $\displaystyle\int_{\mathcal{X}}\mathbb{P}(x)\mathsf{d}\mu>1-\alpha$ and $\mathcal{W}_{J}(\mu,\mu^*)\leq \delta$, where $\mathcal{W}_{J}(\mu,\mu^*)$ is the pseudo-distance (pseudo-metric) between $\mu$ and $\mu^*$ associated with $J$. 
\end{myassump}

Assumption~\ref{assump:P_alpha_opt_sol_seq} is a Slater-type regularity condition in the measure space: it requires the existence of strictly feasible measures arbitrarily close (in $\mathcal{W}_{J}$) to an optimal solution $\mu^*\in A_\alpha$.
It holds, for example, if $\mu^*$ is strictly feasible; more generally, if there exists any strictly feasible $\bar{\mu}$ with $\int_{\mathcal{X}}\mathbb{P}(x)\mathsf{d}\bar{\mu}>1-\alpha$, then $\mu_\varepsilon:=(1-\varepsilon)\mu^*+\varepsilon\bar{\mu}$ is strictly feasible for all $\varepsilon>0$ and satisfies $\mathcal{W}_{J}(\mu_\varepsilon,\mu^*)\to 0$ as $\varepsilon\rightarrow 0$. Assumption~2 may fail in degenerate situations where no strictly feasible measure exists, i.e., when every feasible measure satisfies the chance constraint only at equality. In that case, the feasible set has empty interior in the topology induced by $\mathcal{W}_J$.

Let 
\begin{equation}
    \mathcal{U}_{\mathsf{m}}:=
\Big\{\mu_{\mathsf{m}}\in[0,1]^2:
\sum_{i=1}^2\mu_{\mathsf{m}}^{(i)}=1\Big\}.
\end{equation}
Let $\bm{z}_{\mathsf{m}}:=(\mu_{\mathsf{m}}^{(1)},\mu_{\mathsf{m}}^{(2)},x^{(1)},x^{(2)})$ be a variable in the set $\mathcal{Z}_{\mathsf{m}}:=\mathcal{U}_{\mathsf{m}}\times\mathcal{X}^2$. Then, the optimization problem on $\bm{z}_{\mathsf{m}}$ can be written as
\begin{align} 
    &\min_{\bm{z}_{\mathsf{m}}\in \mathcal{Z}_{\mathsf{m}}} \,\, \sum_{i=1}^2 J(x^{(i)})\mu_{\mathsf{m}}^{(i)}, \tag{$W_{\alpha}$} \label{eq:W_alpha} \\
    &{\normalfont \mathsf{s.t.}}\quad \sum_{i=1}^{2}\mathbb{P}(x^{(i)})\mu_{\mathsf{m}}^{(i)}\geq 1-\alpha. \nonumber
\end{align} 
For Problem \ref{eq:W_alpha}, let us define the feasible set $Z_{\mathsf{m},\alpha}$ of Problem \ref{eq:W_alpha} by
\begin{equation*}
    \mathcal{Z}_{\mathsf{m},\alpha}:=
    \Big\{\bm{z}_{\mathsf{m}}\in\mathcal{Z}_{\mathsf{m}}:\sum_{i=1}^{2}\mathbb{P}(x^{(i)})\mu_{\mathsf{m}}^{(i)}\geq 1-\alpha\Big\}.
\end{equation*}
In addition, define the optimal objective value $\mathcal{J}^{\mathsf{w}}_{\alpha}$ and optimal
solution set $ D_{\alpha}$ of Problem \ref{eq:W_alpha} by
\begin{align*}
    \mathcal{J}^{\mathsf{w}}_{\alpha}
    &:=\min\Big\{\sum_{i=1}^2 J(x^{(i)})\mu_{\mathsf{m}}^{(i)}:\bm{z}_{\mathsf{m}}\in\mathcal{Z}_{\mathsf{m},\alpha}\Big\},\\
    D_{\alpha}&:=\Big\{\bm{z}_{\mathsf{m}}\in\mathcal{Z}_{\mathsf{m},\alpha}:\sum_{i=1}^2 J(x^{(i)})\mu_{\mathsf{m}}^{(i)}=\mathcal{J}^{\mathsf{w}}_{\alpha}\Big\}.
\end{align*}

\begin{myrmk}
    Problem \ref{eq:W_alpha} cannot be regarded as a special case of Affine Chance Constrained Stochastic Programs (ACC-SP) presented in \cite{Cui:2022}, despite both problems involving affine combinations of two violation probabilities. In ACC-SP, the coefficients of the affine combinations are fixed, and the optimization is performed solely over the decision variable $x$. In contrast, Problem \ref{eq:W_alpha} optimizes the probability measures $\mu_{\mathsf{m}}^{(1)}$ and $\mu_{\mathsf{m}}^{(2)}$ together with $x^{(1)}$ and $x^{(2)}$. 
\end{myrmk}

The following theorem shows that the optimal objective values of Problems \ref{eq:W_alpha} and \ref{eq:P_alpha} are equal.

\begin{mytheo}
    \label{theo:P_alpha_2_point_approx}
    Suppose Assumptions \ref{assump:J_h}-\ref{assump:P_alpha_opt_sol_seq} hold. Then, we have $\mathcal{J}^{\mathsf{w}}_{\alpha}=\mathcal{J}^*_{\alpha}$.
\end{mytheo}

The proof of Theorem \ref{theo:P_alpha_2_point_approx} is summarized in Appendix \ref{proof:theo_P_alpha_2_point_approx}. 
To proceed with the proof, some preparations are provided in Appendix \ref{appen:preparation Thm 2}. 
We outline the proof idea here.
First, we show that Problem \ref{eq:P_alpha} can be approximated by a sequence of finite-support probabilistic decisions, and that the corresponding optimal objective values converge to $\mathcal{J}^*_{\alpha}$ (Theorem \ref{theo:sample_approx_convergence}). 
Next, using the preparation results in Appendix \ref{appen:preparation Thm 2}, we prove that the optimal value can always be attained by a probabilistic decision whose probability measure is concentrated on two points. 
This yields the two-point reduction and establishes the equivalence between Problem \ref{eq:P_alpha} and the reduced Problem \ref{eq:W_alpha}. 

For detailed derivations and formal proofs, please refer to Appendices \ref{appen:preparation Thm 2} and \ref{proof:theo_P_alpha_2_point_approx}.

Theorem \ref{theo:P_alpha_2_point_approx} serves as the main result for establishing the reduced problem \ref{eq:W_alpha}, which is equivalent to Problem \ref{eq:P_alpha}. Specifically, instead of solving Problem \ref{eq:P_alpha} in an infinite-dimensional space, we can solve Problem \ref{eq:W_alpha}, whose domain is $[0,1]^2\times\mathcal{X}^2$ in a $(2n+2)$-dimensional space, to obtain optimal probabilistic decisions. 

Let $\bm{z}_{\mathsf{m}}^*:=(\mu_{\mathsf{m}}^{*(1)},\mu_{\mathsf{m}}^{*(2)},x^{*(1)},x^{*(2)})\in D_\alpha$ be an optimal solution of Problem \ref{eq:W_alpha}, where $\mu_{\mathsf{m}}^{*(1)} $ and $ \mu_{\mathsf{m}}^{*(2)}$ are assigned to $ x^{*(1)} $ and $x^{*(2)}$, respectively. Define a mapping  $\Phi$ from the optimal solution set $D_\alpha$ of Problem \ref{eq:W_alpha} to the domain $M(\mathcal{X})$ of Problem \ref{eq:P_alpha} as
$\Phi:D_\alpha\rightarrow M(\mathcal{X})$
as follows: $\forall \bm{z}^*_{\mathsf{m}}\in D_\alpha$, assign a probability measure in $ M(\mathcal{X})$, 
$\tilde{\mu}_{\bm{z}^*_{\mathsf{m}}}=\Phi(\bm{z}^*_{\mathsf{m}})\in M(\mathcal{X})$,
with $\tilde{\mu}_{\bm{z}^*_{\mathsf{m}}}(\{x^{*(1)}\})=\mu_{\mathsf{m}}^{*(1)}$ and $\tilde{\mu}_{\bm{z}^*_{\mathsf{m}}}(\{x^{*(2)}\})=\mu_{\mathsf{m}}^{*(2)}$. 

Then, by Theorem \ref{theo:P_alpha_2_point_approx}, we have
\begin{equation}
    \label{eq:tilde_mu_star_m_J}
    \int_{\mathcal{X}}J(x)\mathsf{d}\tilde{\mu}_{\bm{z}^*_{\mathsf{m}}}(x)=\sum_{i=1}^2J(x^{*(i)})\mu_{\mathsf{m}}^{*(i)}=\mathcal{J}^*_\alpha.
\end{equation}
In addition, if $\bm{z}_{\mathsf{m}}^*\in D_\alpha$ is the optimal solution of Problem \ref{eq:W_alpha}, we have
\begin{equation}
    \label{eq:tilde_mu_star_m_alpha}
    \int_{\mathcal{X}}\mathbb{P}(x)\mathsf{d}\tilde{\mu}_{\bm{z}^*_{\mathsf{m}}}(x)=\sum_{i=1}^2 \mathbb{P}(x^{*(i)})\mu_{\mathsf{m}}^{*(i)} \geq 1-\alpha.
\end{equation}

From \eqref{eq:tilde_mu_star_m_J} and \eqref{eq:tilde_mu_star_m_alpha}, we have that $\tilde{\mu}_{\bm{z}^*_{\mathsf{m}}}$ is within the optimal set $A_\alpha$ of Problem \ref{eq:P_alpha}. Specifically, by solving Problem~\ref{eq:W_alpha}, we can obtain the optimal solution $\tilde{\mu}_{\bm{z}^*_{\mathsf{m}}}$ of Problem \ref{eq:P_alpha} through the mapping $\Phi$. 

We next summarize the properties of the optimal solutions of Problem \ref{eq:W_alpha}. 

\begin{mycoro}
    \label{coro:z_m_opt}
    Let $\bm{z}_{\mathsf{m}}^*=(\mu_{\mathsf{m}}^{*(1)},\mu_{\mathsf{m}}^{*(2)},x^{*(1)},x^{*(2)})\in D_\alpha$ be an  optimal solution of Problem \ref{eq:W_alpha}. Let $\tilde{\alpha}^{(1)}=1-\mathbb{P}(x^{*(1)})$ and $\tilde{\alpha}^{(2)}=1-\mathbb{P}(x^{*(2)})$. Then, $x^{*(1)}\in X_{\tilde{\alpha}^{(1)}} $ and $ x^{*(2)}\in X_{\tilde{\alpha}^{(2)}}$. 
    Recall that for any violation level $\bar{\alpha}\in[0,1]$, $\mathcal{X}_{\bar{\alpha}}$ is the feasible set of Problem \ref{eq:Q_alpha} at level $\bar{\alpha}$ (i.e., the set defined in \eqref{eq:mathX_alpha_Q}).
\end{mycoro}

The proof of Corollary \ref{coro:z_m_opt} is presented in Appendix \ref{proof:coro_z_m_opt}. By Corollary \ref{coro:z_m_opt}, the optimal solution of Problem \ref{eq:W_alpha} contains two points in $\mathcal{X}$ such that are optimal solutions of $Q_{\tilde{\alpha}^{(1)}}$ and $Q_{\tilde{\alpha}^{(2)}}$, respectively. The following corollary summarizes the connections between Problems \ref{eq:P_alpha} and~\ref{eq:Q_alpha}.

\begin{mycoro}
    \label{coro:W_alpha_Q_alpha}
    If $J^*_{\tilde{\alpha}}$ is a convex function of $\tilde{\alpha}$ or when $\alpha=0$, then $\mathcal{J}^*_{\alpha}=J^*_\alpha$.
\end{mycoro}

The proof of Corollary \ref{coro:W_alpha_Q_alpha} is provided in Appendix \ref{proof:coro_W_alpha_Q_alpha}. Corollary \ref{coro:W_alpha_Q_alpha} shows that it is not necessary to consider CCPMO for improving the expected performance in two cases: (1) when $J^*_{\tilde{\alpha}}$ is a convex function of $\tilde{\alpha}$, or (2) when $\alpha=0$. However, verifying whether the function $J^*_{\tilde{\alpha}}$ is convex in $\tilde{\alpha}$ a prior is generally not feasible. On the other hand, when $\alpha=0$, the constraint $\bar{h}(x,\xi)$ should be satisfied with probability 1 (w. p. 1). In this case, the deterministic decision can achieve the optimal solution, and the probabilistic decision is unnecessary. 

\section{Sample-based smooth approximation for solving Problem \ref{eq:W_alpha}}
\label{sec:approximation}

Section \ref{sec:reduction} has demonstrated that Problem \ref{eq:P_alpha} admits an optimal probabilistic decision, where the probability measure is concentrated on two points. By solving Problem \ref{eq:W_alpha} in a $(2n+2)$-dimensional space, we can obtain an optimal solution of Problem \ref{eq:P_alpha}. However, due to the presence of chance constraints induced by the random variable $\xi$ in \eqref{eq:def_P_x} and \eqref{eq:def_bar_h}, Problem \ref{eq:W_alpha} remains computationally intractable. To address this issue, we propose an approximate method for solving Problem \ref{eq:W_alpha} by leveraging sample extraction from the support set $\Xi$. 

\subsection{Uniform convergence of optimal objective values}

Refer to \cite{Pena:2020}, we have the following assumption on Problem \ref{eq:Q_alpha} for every $\alpha\in(0,1]$.

\begin{myassump}\label{assump:opti_solu_Q_alpha}
    For any $\alpha\in (0,1]$, there exists a globally optimal solution of Problem \ref{eq:Q_alpha}, $\bar{x}_\alpha$, such that for any $\eta>0$ there is an $x\in\mathcal{X}$ satisfying $\|x-\bar{x}_\alpha\|\leq\eta$ and $\mathbb{P}(x)>1-\alpha$. 
\end{myassump}

Assumption \ref{assump:opti_solu_Q_alpha} implies that there exists a sequence $\{x_k\}_{k=1}^\infty\subset\mathcal{X}$ that converges to an optimal solution $\bar{x}_\alpha$ 
of Problem \ref{eq:Q_alpha}, where each $x_k$ satisfies $\mathbb{P}(x_k)>1-\alpha$. 

Recall that the constraint function $\bar{h}(x,\xi)$ in \eqref{eq:def_bar_h} depends on the random variable $\xi$, whose sample space is denoted by $\Xi$. Let $\mathcal{D}_N=\{\xi^{(1)},\ldots,\xi^{(N)}\}$ be a set of $N$ samples, where $\xi^{(j)},j=1,\ldots,N$ is drawn from $\Xi$, the support of $\xi$. With dataset $\mathcal{D}_N$, we can define a sample-based approximation of $\mathbb{P}(x)$ by
\begin{equation}
    \label{eq:P_N_x}
    \hat{\mathbb{P}}_N(x):=\frac{1}{N}\sum_{j=1}^N\mathbb{I}\left(\bar{h}(x,\xi^{(j)})\right).
\end{equation}
One limitation of the sample-based approximation approach is its lack of differentiability, which arises from the discontinuous nature of the indicator function $\mathbb{I}(\cdot)$. To overcome this, as in \cite{Pena:2020}, for a positive parameter $\varepsilon$, a function $\Lambda_{\varepsilon}(\cdot)$ is defined
\begin{equation*}
    \Lambda_{\varepsilon}(y)=
     \left\{
     \begin{array}{ll}
          0,  & \mbox{if}\ y\geq \varepsilon,\\
          \lambda_{\varepsilon}(y), & \mbox{if}\ -\varepsilon<y<\varepsilon, \\
          1, & \mbox{otherwise},
          \end{array}
     \right.
\end{equation*}
where $\lambda_{\varepsilon}:[-\varepsilon,\varepsilon]\rightarrow[0,1]$ is a symmetric and strictly decreasing smooth function that makes $\Lambda_{\varepsilon}(\cdot)$ a continuously differentiable function. 

For example, one may choose
\begin{equation*}
\Lambda_{\epsilon}(y)=
\begin{cases}
0, & y\ge \epsilon,\\[2mm]
\frac{1}{2}-\frac{3}{4}\frac{y}{\epsilon}+\frac{1}{4}\left(\frac{y}{\epsilon}\right)^3, & -\epsilon<y<\epsilon,\\[2mm]
1, & y\le -\epsilon,
\end{cases}
\end{equation*}
which is symmetric in the sense that $\Lambda_{\epsilon}(-y)=1-\Lambda_{\epsilon}(y)$ and is continuously differentiable on $\mathbb{R}$.

With a smooth parameter $\varepsilon>0$, sample number $N>0$, a positive parameter $\gamma>0$, we define $\tilde{\mathbb{P}}^{\gamma,\varepsilon}_N(x)$ by
\begin{equation}
    \tilde{\mathbb{P}}^{\gamma,\varepsilon}_N(x):=\frac{1}{N}\sum_{j=1}^N\Lambda_{\varepsilon}\left(\bar{h}(x,\xi^{(j)})+\gamma\right).
\end{equation}
The function $\tilde{\mathbb{P}}^{\gamma,\varepsilon}_N(x)$ is a sample-based smooth approximation of $\mathbb{P}(x)$. Compared to the sample-based approximation $\hat{\mathbb{P}}_{N}(x)$, $\tilde{\mathbb{P}}^{\gamma,\varepsilon}_N(x)$ is differentiable and gives a smooth approximation of the boundary of the feasible region.

With a given set $\mathcal{D}_N$, positive parameters $\alpha'\in[0,\alpha),\gamma,\varepsilon>0$, a sample-based smooth approximate problem of \ref{eq:W_alpha} is formulated as:
\begin{align} \label{eq:til_W_alpha_gamma_N}
    &\min_{\bm{z}_{\mathsf{m}} \in \mathcal{Z}_{\mathsf{m}}} \sum_{i=1}^2 J(x^{(i)})\mu_{\mathsf{m}}^{(i)}, \tag{$\tilde{W}_{\alpha'}^{\gamma,\varepsilon}(\mathcal{D}_N)$} \\
    &{\normalfont \mathsf{s.t.}}\quad \sum_{i=1}^{2}\tilde{\mathbb{P}}^{\gamma,\varepsilon}_N(x^{(i)})\mu_{\mathsf{m}}^{(i)}\geq 1-\alpha'. \nonumber
\end{align} 
Let $$\mathcal{J}_{\mathsf{m}}(\bm{z}_{\mathsf{m}}):=\displaystyle\sum_{i=1}^2 J(x^{(i)})\mu_{\mathsf{m}}^{(i)},\ \mathcal{P}_{\mathsf{m}}(\bm{z}_{\mathsf{m}}) :=\displaystyle\sum_{i=1}^{2}\tilde{\mathbb{P}}^{\gamma,\varepsilon}_N(x^{(i)})\mu_{\mathsf{m}}^{(i)}.$$
The feasible region $\mathcal{Z}_{\mathsf{m},\alpha'}^{\gamma,\varepsilon}(\mathcal{D}_N)$, the optimal objective value $\tilde{\mathcal{J}}_{\alpha'}^{\gamma,\varepsilon}(\mathcal{D}_N)$, and the optimal solution set $\tilde{D}_{\alpha'}^{\gamma,\varepsilon}(\mathcal{D}_N)$ of Problem \ref{eq:til_W_alpha_gamma_N} are defined by 
\begin{align*}
    \mathcal{Z}_{\mathsf{m},\alpha'}^{\gamma,\varepsilon}(\mathcal{D}_N)&:=
    \Big\{\bm{z}_{\mathsf{m}}\in \mathcal{Z}_{\mathsf{m}}:\mathcal{P}_{\mathsf{m}}(\bm{z}_{\mathsf{m}})\geq 1-\alpha'\Big\},  \\
    \tilde{\mathcal{J}}_{\alpha'}^{\gamma,\varepsilon}(\mathcal{D}_N)&:=\min_{\bm{z}_{\mathsf{m}}\in\mathcal{Z}_{\mathsf{m},\alpha'}^{\gamma,\varepsilon}(\mathcal{D}_N)}\mathcal{J}_{\mathsf{m}}(\bm{z}_{\mathsf{m}}), \\
   \tilde{D}_{\alpha'}^{\gamma,\varepsilon}(\mathcal{D}_N)&:=
    \!\left\{
    \bm{z}_{\mathsf{m}}\in\mathcal{Z}_{\mathsf{m},\alpha'}^{\gamma,\varepsilon}(\mathcal{D}_N):\mathcal{J}_{\mathsf{m}}(\bm{z}_{\mathsf{m}})=\tilde{\mathcal{J}}_{\alpha'}^{\gamma,\varepsilon}(\mathcal{D}_N)\right\}. 
\end{align*}
The uniform convergence in the sense of almost everywhere on $\tilde{\mathcal{J}}_{\alpha'}^{\gamma,\varepsilon}(\mathcal{D}_N)$ and $\tilde{D}_{\alpha'}^{\gamma,\varepsilon}(\mathcal{D}_N)$ is summarized in the following theorem. 
\begin{mytheo}\label{theo:uniform_convergence_W}
    Suppose Assumptions \ref{assump:J_h}-\ref{assump:opti_solu_Q_alpha} hold. If $\gamma=0$, then, as $N\rightarrow\infty$ and $\varepsilon\rightarrow 0$, $\tilde{\mathcal{J}}_{\alpha}^{\gamma,\varepsilon}(\mathcal{D}_N)\rightarrow\mathcal{J}^*_{\alpha}$ and $\mathbb{D}\left(\tilde{D}_{\alpha}^{\gamma,\varepsilon}(\mathcal{D}_N),D_\alpha\right)\rightarrow 0$~ w. p. 1,  where
     $$\mathbb{D}(Z_1,Z_2):=\sup_{z\in Z_1}\inf_{z'\in Z_2}\|z-z'\|$$ is the deviation of set $Z_1$ from set $Z_2$.
\end{mytheo}

The proof of Theorem \ref{theo:uniform_convergence_W} is presented in Appendix \ref{proof:theo_uniform_convergence_W}. Theorem \ref{theo:uniform_convergence_W} demonstrates that the objective function $\tilde{\mathcal{J}}_{\alpha}^{\gamma,\varepsilon}(\mathcal{D}_N)$ of the approximate problem \ref{eq:til_W_alpha_gamma_N} converges to the optimal cost $\mathcal{J}^*_{\alpha}$ of the original problem \ref{eq:P_alpha} w. p. 1 as the sample size $N$ goes to infinity and the smooth parameter $\varepsilon$ decreases to zero. In addition, the optimal solution set $\tilde{D}_{\alpha}^{\gamma,\varepsilon}(\mathcal{D}_N)$ of the approximate problem \ref{eq:til_W_alpha_gamma_N} converges to a subset of the optimal solution set $D_\alpha$ of the original problem \ref{eq:P_alpha} w. p. 1 as the sample size $N$ goes to infinity and the smooth parameter $\varepsilon$ decreases to zero.

\begin{myrmk}\label{remark:theorem 3}
    Problem \ref{eq:til_W_alpha_gamma_N} is a non-convex problem, and existing algorithms may lead to local optima rather than global optima. The uniform convergence of optimality in Theorem \ref{theo:uniform_convergence_W} pertains to global optimality. However, this convergence result can also be extended to local optimality case as we can consider a subset of $\mathcal{X}$ to make local optimal solution become global optima in it.
\end{myrmk}

\subsection{Feasibility with finite samples}

We next discuss the condition under which an optimal solution for Problem \ref{eq:til_W_alpha_gamma_N} is feasible for Problem \ref{eq:W_alpha}. Define $\hat{\mathbb{P}}^{\gamma,\varepsilon}(x)$ for $x\in\mathcal{X}$ by
\begin{equation}
    \hat{\mathbb{P}}^{\gamma,\varepsilon}(x):=\int_{\Xi}\Lambda_{\varepsilon}\left(\bar{h}(x,\xi)+\gamma\right)p(\xi)\mathsf{d}\xi.
\end{equation}
We have the probabilistic feasibility guarantee as stated in the following theorem.
\begin{mytheo}
    \label{theo:prob_fea_gua}
     Let $\tilde{\bm{z}}_{\mathsf{m},\alpha'}\in\tilde{D}^{\gamma,\varepsilon}_{\alpha'}(\mathcal{D}_N)$ be an optimal solution of Problem \ref{eq:til_W_alpha_gamma_N}, which is specified as $\tilde{\bm{z}}_{\mathsf{m},\alpha'}:=(\tilde{\mu}_{\mathsf{m}}^{(1)},\tilde{\mu}_{\mathsf{m}}^{(2)},\tilde{x}^{(1)},\tilde{x}^{(2)})$. Then, if $\alpha'<\alpha$, we have
    \begin{equation}
        \mathsf{Pr}\{\tilde{\bm{z}}_{\mathsf{m},\alpha'}\notin \mathcal{Z}_{\mathsf{m},\alpha}\}\leq\exp\{-2N\left(R^{\gamma,\varepsilon}_{\alpha'}\right)^2\},
    \end{equation}    
    where
    \begin{equation}\label{eq:R_def}
        R^{\gamma,\varepsilon}_{\alpha'}:= \!\inf_{\bm{z}_{\mathsf{m}}\in\mathcal{Z}_{\mathsf{m}}} \sum_{i=1}^2\left(\mathbb{P}(x^{(i)})-\hat{\mathbb{P}}^{\gamma,\varepsilon}(x^{(i)})\right)\mu_{\mathsf{m}}^{(i)}+(\alpha-\alpha').
    \end{equation}
\end{mytheo}

The proof of Theorem \ref{theo:prob_fea_gua} is presented in Appendix \ref{proof:theo_prob_fea_gua}. Theorem \ref{theo:prob_fea_gua} shows that if the violation probability threshold of the approximate Problem \ref{eq:til_W_alpha_gamma_N} is appropriately reduced to $\alpha'<\alpha$, then an optimal solution $\tilde{\bm{z}}_{\mathsf{m},\alpha'}$ of the approximate problem \ref{eq:til_W_alpha_gamma_N} is not a feasible solution for Problem \ref{eq:W_alpha} with a probability that decays exponentially as the sample size $N$ increases.

\begin{myrmk}
    As discussed in Remark \ref{remark:theorem 3}, the probabilistic feasibility result in Theorem \ref{theo:prob_fea_gua} also holds with local optimal solution of Problem \ref{eq:til_W_alpha_gamma_N}.
\end{myrmk}

\subsection{Proposed algorithms}

Theorem \ref{theo:P_alpha_2_point_approx} demonstrates that an optimal solution of Problem \ref{eq:P_alpha} can be obtained by solving Problem \ref{eq:W_alpha}, as Problem \ref{eq:P_alpha} admits an optimal probabilistic decision where the probability measure is concentrated on two points. Furthermore, Theorem \ref{theo:uniform_convergence_W} establishes that the solution of a sample-based approximate problem \ref{eq:til_W_alpha_gamma_N} converges to an optimal solution of Problem \ref{eq:W_alpha} as the sample size $N$ the random variable $\xi$ goes to infinity. Additionally, Theorem \ref{theo:prob_fea_gua} guarantees that the optimal solution of Problem \ref{eq:til_W_alpha_gamma_N} with a finite sample size $N$ satisfies the original chance constraint of Problem \ref{eq:W_alpha} with a probability bound. 

Leveraging these theoretical results, Theorems \ref{theo:P_alpha_2_point_approx}, \ref{theo:uniform_convergence_W}, and \ref{theo:prob_fea_gua}, the algorithm to obtain an approximate solution of Problem \ref{eq:P_alpha} is proposed in Algorithm \ref{alg:random_opt_tilde}. Furthermore, the procedure for implementing a probabilistic decision in a repeated decision process based on a solution $\tilde{\mu}_{\mathsf{a}}$ is summarized in Algorithm \ref{alg:prob_decision}.

\begin{algorithm}[h]
    \caption{Proposed algorithm for obtaining an approximate solution of Problem \ref{eq:P_alpha}} \label{alg:random_opt_tilde}
    \begin{algorithmic}[1]
        \Require  Parameter $\gamma>0$; sample number $N$ of $\xi$; smooth degree $\varepsilon>0$; probability level $\alpha'$
        \Ensure Approximate solution $\mu^*_{\alpha}\in M(\mathcal{X})$ that satisfies $\mu^*_{\alpha}(\{\tilde{x}_\mathsf{m}^{(1)}\})=\tilde{\mu}_\mathsf{m}^{(1)}$ and $\mu^*_{\alpha}(\{\tilde{x}_\mathsf{m}^{(2)}\})=\tilde{\mu}_\mathsf{m}^{(2)}$
        \Statex 
        \State Randomly extract $\mathcal{D}_N$ from $\Xi$ 
        \State Solve Problem \ref{eq:til_W_alpha_gamma_N} and obtain $\tilde{\bm{z}}_{\mathsf{m}}=(\tilde{\mu}_\mathsf{m}^{(1)},\tilde{\mu}_\mathsf{m}^{(2)},\tilde{x}_\mathsf{m}^{(1)},\tilde{x}_\mathsf{m}^{(2)})$
    \end{algorithmic}
\end{algorithm}

\begin{algorithm}[h]
    \caption{Implementation and evaluation of $\mu^*_{\alpha}$-based probabilistic decision}
    \label{alg:prob_decision}
    \begin{algorithmic}[1]
      \Require Solution $\mu^*_{\alpha}$; repeat time of decision process~$T$
      \Statex 
      \For{$t = 1, \dots, T$} 
        \State Randomly generate $\kappa$ from $[0,1]$ according to uniform distribution
        \If{$\kappa \leq \mu^*_{\alpha}(\{\tilde{x}_\mathsf{m}^{(1)}\})$}
          \State Implement $x(t) = \tilde{x}^{(1)}_{\mathsf{m}}$
        \ElsIf{$\kappa > \mu^*_{\alpha}(\{\tilde{x}_\mathsf{m}^{(1)}\})$}
          \State Implement $x(t) = \tilde{x}^{(2)}_{\mathsf{m}}$ 
        \EndIf
        \State Calculate instant cost $J(x(t))$
        \State Calculate instant risk indicator $\mathbb{I}(\bar{h}(x(t),\xi(t)))$ 
      \EndFor
      \State Compute total cost: $\displaystyle\sum_{t=1}^T J(x(t))$
      \State Compute average risk: $\frac{1}{T} \displaystyle\sum_{t=1}^T \mathbb{I}(\bar{h}(x(t),\xi(t)))$
    \end{algorithmic}
\end{algorithm}

\section{Numerical Example}\label{sec:validations}

Numerical validations have been implemented to compare the performance of the proposed probabilistic decision-making method and several existing deterministic decision-making methods. The chosen application case study is a chance-constrained trajectory planning problem with obstacles, a widely recognized challenge in trajectory planning for autonomous driving.  

\subsection{Model and simulation settings}

Consider a quadrotor system control problem in turbulent conditions. The control problem is expressed as follows:
\begin{align} \label{eq:P_QSC}
    &\underset{\mu\in M(\mathcal{U}_T)}{{\normalfont\mathsf{min}}} \,\, \mathcal{J}(\mu)=\mathbb{E}\left\{\ell^s(s)+\ell^u(u)\right\}, \tag{$P_{\mathsf{QSC}}$}\\
    &{\normalfont \mathsf{s.t.}}\quad  s_{t+1}=As_t+B(m)u_t+d(s_t,\varphi)+\omega_t, \nonumber\\
    & \ \ \ \ \quad u\sim M(\mathcal{U}_T),\; t=0,1,...,T-1,\nonumber\\
    &\ \ \ \ \quad \mathsf{Pr}\{\left(\land_{t=1}^{T-1}s_t\notin\mathcal{O}\right)\land(s_T\in\mathcal{S}_{\mathsf{goal}})\} \geq 1-\alpha,\nonumber
\end{align}
where the cost functions are given by
\begin{align*}
    \ell^s(s) &= \frac{1}{T}\sum_{t=0}^{T-1}(p_{x,t+1}-p_{x,t})^2+(p_{y,t+1}-p_{y,t})^2,\\
    \ell^u(u) &= \frac{0.1}{T}\sum_{t=0}^{T-1}u_{1,t}^2+u_{2,t}^2.
\end{align*}

\begin{figure*}
    \centering
    \includegraphics[width=\hsize]{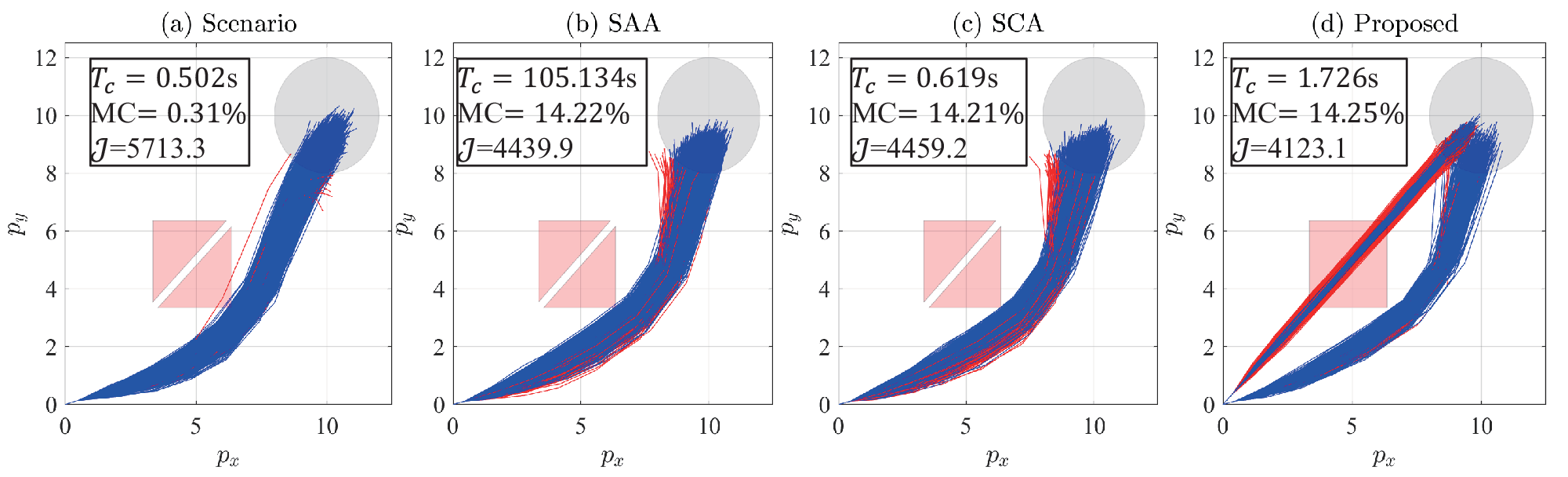}
    \centering
    \caption{Solutions from different methods for the tolerable failure probability threshold $\alpha=15\%$. Blue trajectories from Monte-Carlo (MC) simulations denote feasible trajectories that reach the goal set $\mathcal{S}_{goal}$ and avoid obstacles $\mathcal{O}$. Red trajectories mean that either constraint is violated in the simulation. The expression $\mathsf{MC}=14.22\%$ means the probability of having red trajectories in the MC simulations. (a) \textbf{Scenario}; (b) \textbf{SAA}; (c) \textbf{SCA}; (d) \textbf{Proposed}.}
    \label{fig:control_exa}
\end{figure*}

The parameters $A$, $B(m)$, $d(s_t,\varphi)$ are given by
\begin{equation*}
    A=
    \begin{bmatrix}
        1 & \Delta t & 0 & 0 \\
        0 & 1 & 0 & 0 \\
        0 & 0 & 1 & \Delta t \\
        0 & 0 & 0 & 1
    \end{bmatrix},\;
    B(m)=\frac{1}{m}
    \begin{bmatrix}
        \frac{\Delta t^2}{2} & 0 \\
        \Delta t & 0 \\
        0 & \frac{\Delta t^2}{2} \\
        0 & \Delta t
    \end{bmatrix},
\end{equation*}
and $d(s_t,\varphi)=-\varphi[\frac{\Delta t^2|v_x|v_x}{2}\ \Delta t|v_x|v_x\ \frac{\Delta t^2|v_y|v_y}{2}\ \Delta t|v_y|v_y]^{\top}$. $\Delta t=1$ is the sampling time. The system state is denoted as $s_t=[p_{x,t},v_{x,t},p_{y,t},v_{y,t}]^{\top} \in\mathbb{R}^4$ and the control input is denoted as $u_t=[u_x,u_y]^{\top} \in \mathbb{R}^2 $. Then, the state and control trajectories are denoted as $s=\{s_t\}_{t=1}^{T}$ and $u=\{u_t\}_{t=1}^{T-1}$. The system starts from an initial point $s_0=[-0.5,0,-0.5,0]$. 
The system is expected to reach the destination set $\mathcal{S}_{\mathsf{goal}}=\{(p_x,v_x,p_y,v_y)|\|(p_x-10,p_y-10)\|\leq 2, v_x=0,v_y=0\}$ at time $t=T=10$ while avoiding two polytopic obstacles $\mathcal{O}$ shown in Fig. \ref{fig:control_exa}. 
The dynamics are parametrized by uncertain parameter vector $\xi_t=[m,\varphi]^\top$, where $m>0$ represents the system's mass and $\varphi>0$ is an uncertain drag coefficient. 
The parameter vector $\xi$ of the system is uncorrelated random variables such that $(m-0.75)/0.5\sim \mathsf{Beta}(2,2)$ and $(\varphi-0.4)/0.2\sim \mathsf{Beta}(2,5)$, where $\mathsf{Beta}(a,b)$ denotes a Beta distribution with shape parameters $(a,b)$. $\omega_t\in\mathbb{R}^4$ is the uncertain disturbance at time step $t$, which obeys multivariate normal distribution with zero means and a diagonal covariance matrix $\Sigma$ with the diagonal elements as $\Sigma(1,1)=0.01,\ \Sigma(2,2)=0.75,\ \Sigma(3,3)=0.01,\ \Sigma(4,4)=0.75$.

By solving Problem \ref{eq:P_QSC}, we can get probability measure $\mu$ to implement a probabilistic decision. To obtain deterministic decision for this example, we can solve the following problem instead.
\begin{align} 
\label{eq:P_QSC_det}
&\underset{u\in \mathcal{U}_T}{{\normalfont\mathsf{min}}} \,\, J(u)=\mathbb{E}\{\ell^s(s)+\ell^u(u)\}, \tag{$P_{\mathsf{QSC}}^{\mathsf{det}}$}\\
&{\normalfont \mathsf{s.t.}}\quad  s_{t+1}=As_t+B(m)u_t+d(s_t,\varphi)+\omega_t,  \nonumber\\
& \ \ \ \ \qquad t=0,1,...,T-1,\nonumber\\
&\ \ \ \ \quad \mathsf{Pr}\{\left(\land_{t=1}^{T-1}s_t\notin\mathcal{O}\right)\land(s_T\in\mathcal{S}_{\mathsf{goal}})\} \geq 1-\alpha.\nonumber
\end{align}

\noindent In Problem \ref{eq:P_QSC_det}, the optimization variable is the deterministic open-loop input sequence $u\in\mathcal{U}_T$. The expectation in the objective is taken with respect to the exogenous uncertainties in the system dynamics, i.e., the uncertain parameters $\xi_t=[m,\varphi]^\top$ and the disturbance $\omega_t$, under their prescribed distributions, while $u$ is fixed deterministically during optimization.

With this example, we would like to compare the performance of the following methods: 
\begin{itemize}
    \item \textbf{Proposed}: Solve \ref{eq:P_QSC} and implement probabilistic decisions by Algorithms \ref{alg:random_opt_tilde} and \ref{alg:prob_decision}; 
    \item \textbf{Scenario}: Solve \ref{eq:P_QSC_det} by scenario approach \cite{Campi_unconvex} and implement deterministic decision; 
    \item \textbf{SAA}: Solve \ref{eq:P_QSC_det} by sample average approach \cite{Luedtke:2008} and implement deterministic decision; 
    \item \textbf{SCA}: Solve \ref{eq:P_QSC_det} by sample-based smooth approximation \cite{Pena:2020} and implement deterministic decision.
\end{itemize}

\subsection{Results and discussions}

The simulation results are shown in Fig. \ref{fig:control_exa}, comparing the above four methods under the same conditions: the violation probability threshold $\alpha= 15\% $, the number of samples $N = 2000$. Additionally, the parameters $\varepsilon$ and $\gamma$ are both set as $0.01$ for \textbf{SCA} and \textbf{Proposed}. Fig. \ref{fig:control_exa} shows the outcomes of $10,000$ Monte-Carlo (MC) simulations of the quadrotor system using the open-loop policies computed by four methods. In this work, we employ the interior-point algorithm presented in \cite{Vanderbei:1999} to solve \textbf{Proposed}.

The deterministic decision generated by \textbf{Scenario} yields a highly safe trajectory with only $0.31\%$ constraint violation in the MC simulations. However, \textbf{Scenario} results in a high expected cost $\mathcal{J}$. This conservatism arises because \textbf{Scenario} was inherently designed for robust optimization and thus it often produces conservative solutions that exceed the required violation probability threshold. In contrast, the deterministic control policies generated by \textbf{SAA} and \textbf{SCA} produce similar trajectories, which substantially optimize the cost $\mathcal{J}$ compared to \textbf{Scenario} while still satisfying the required chance constraint. Notably, \textbf{SAA} requires solving a mixed-integer program, which leads to a considerably longer computation time $T_c$ than \textbf{SCA} and \textbf{Scenario}. For the probabilistic decision generated by \textbf{Proposed}, the cost $\mathcal{J}$ is further improved, albeit with a moderate increase in computation time $T_c$, while maintaining the chance constraint to ensure safety. Therefore, this demonstrates the potential of probabilistic control policies to enhance the expected performance with guaranteeing the safety chance constraints. 

Furthermore, a comprehensive statistical analysis was conducted by performing 1,000 MC trials, each using the same parameters as before. In each MC trial, $10,000$ MC simulations of the quadrotor system are executed to compute the violation probability and the expected cost. The computation results are reported in Table \ref{tb:accuracy}.
The deterministic policy obtained by \textbf{Scenario} shows high robustness against uncertainties but incurs higher expected cost. This deterministic decision corresponds to a completely robust decision, which cannot be improved if complete robustness is required. 
However, in our setting, a risk of $ \alpha = 15\% $ violation probability is allowed, which can enable the possibility of better performance using probability measure decisions. 
\textbf{SAA} and \textbf{SCA} achieve better performance on $\mathbb{E}\{\mathcal{J}\}$, while ensuring zero probability of violating the required $15\%$ risk level when the sample size $N$ is at least $1,500$. The deterministic decisions obtained by \textbf{SAA} and \textbf{SCA} converge to the optimal ones under the $\alpha = 15\%$ risk probability. Notably, compared to \textbf{SAA}, \textbf{SCA} achieves reduced computation time since it employs a smooth indicator function to avoid sticking at certain points due to the non-smoothness of the constraint function. 
Compared to the deterministic decisions obtained by \textbf{SAA} and \textbf{SCA}, the probabilistic decision obtained by \textbf{Proposed} further improves the expected objective value $\mathbb{E}\{\mathcal{J}\}$, validating our theoretical results (Theorem \ref{theo:P_alpha_2_point_approx}). 

\begin{table}[t]
 \centering
 \caption{Statistical analysis of the performance. Results (mean) of 1000 trials are reported.}  
 \begin{tabular}{ccccc}
 \hline
 Methods & $N$ & $\mathbb{E}\{\mathcal{J}\}$ & $\mathbb{\mathsf{Pr}}\{\alpha>15\%\}$ & $\mathbb{E}\{T_c\}$ [s]   \\
 \hline
  & $100$ & 4763.1 & 0.023 & \textbf{0.048} \\
 & $200$ & 4928.1 & 0.0042 & \textbf{0.079} \\
 & $500$ & 5189.1 & 0.0012 & \textbf{0.128} \\
 \textbf{Scenario}& $800$  & 5281.7  & 0 & \textbf{0.176}\\
 & $1000$ & 5402.3 & 0 & \textbf{0.288} \\
 & $1500$ & 5562.9 & 0 & \textbf{0.378} \\
 & $2000$ & 5693.5 & 0 & \textbf{0.496} \\
 \hline
  & $100$ & 4088.2 & 0.0394 & 9.76 \\
 & $200$ & 4176.9 & 0.0129 & 14.67 \\
  & $500$ & 4257.4 & 0.0029 & 23.62 \\
 \textbf{SAA}& $800$ & 4319.8  & 0.0023 & 34.22\\
 & $1000$ & 4352.6 & 0.0016 & 49.38 \\
& $1500$ & 4398.3 & 0 & 68.92 \\
 & $2000$ & 4432.1 & 0 & 97.84 \\
 \hline
   & $100$ & 4092.3 & 0.0379 & 0.063 \\
 & $200$ & 4198.7 & 0.0125 & 0.085 \\
  & $500$ & 4268.5 & 0.0018 & 0.142 \\
 \textbf{SCA}& $800$ & 4342.9  & 0.0015 & 0.213\\
 & $1000$ & 4381.1 & 0.0011 & 0.327 \\
& $1500$ & 4429.2 & 0 & 0.479 \\
 & $2000$ & 4465.3 & 0 & 0.631 \\
 \hline
   & $100$ & \textbf{3785.1} & 0.0387 & 0.153 \\
 & $200$ & \textbf{3895.5} & 0.0128 & 0.221 \\
  & $500$ & \textbf{3979.1} & 0.0022 & 0.365 \\
 \textbf{Proposed}& $800$ & \textbf{4025.9}  & 0.0016 & 0.578\\
 & $1000$ & \textbf{4089.7} & 0.0012 & 0.867  \\
& $1500$ & \textbf{4106.1} & 0 & 1.151 \\
 & $2000$ & \textbf{4120.8} & 0 & 1.720 \\
 \hline
 \end{tabular}
 \label{tb:accuracy} 
\end{table}

\section{Conclusions and Future Work}\label{sec:conclusion}

In a probabilistic decision-making problem under chance constraint, decision variables are probabilistically chosen to optimize the expected objective value under chance constraint. We formulate the probabilistic decision-making problem under chance constraints as the chance-constrained probability measure optimization. We prove the existence of the optimal solution to the chance-constrained probability measure optimization. Furthermore, we show that there exists an optimal solution that the probability measure is concentrated on two decisions. Then, an equivalent reduced problem of the chance-constrained probability measure optimization is established. A sample-based smooth approximation method has been proposed to solve the reduced problem. The analysis of uniform convergence and feasibility is given. A numerical example of a quadrotor system control problem has been implemented to validate the performance of the proposed probabilistic decision-making under chance constraints. Future work will be focused on implementing the proposed method to design optimal feedback probabilistic control policy with probabilistic safety guarantees.

\section*{Acknowledgments}

The authors would like to thank Prof. Xiaoping Xue from Department of Mathematics at Harbin Institute of Technology for helping to organize the proof of Theorem \ref{theo:P_alpha_opt_solution}, and discussing proof of Theorems \ref{theo:sample_approx_convergence} and \ref{theo:P_alpha_2_point_approx}. This research was partially supported by JST Moonshot R\&D Grant Number JPMJMS2021. Dr Ye Wang acknowledges support from the Australian Research Council through the Discovery Early Career Researcher Award (DE220100609).

\bibliographystyle{apacite}
\bibliography{autosam}           

\appendix

\section{Preparation for Proof of Theorem \ref{theo:P_alpha_2_point_approx}}\label{appen:preparation Thm 2}

Let $\chi_S=\left(x^{(1)},...,x^{(S)}\right)$ be a multi-sample realization of the augmented space $\mathcal{X}^S$, where $S\in\mathbb{N}$. For any $\chi_S$, we can define a set of discrete probability measures by
\begin{equation*}
    \mathcal{U}_S:=\Big\{\mu_S\in [0,1]^S:\sum_{i=1}^S\mu_S^{(i)}=1\Big\}.
\end{equation*}
The realization $\chi_S$ forms a finite sample space when equipped with a discrete probability measure $\mu_S\in\mathcal{U}_S$. Here the $i$-th element $\mu_S^{(i)}$ denotes the probability of selecting decision $x^{(i)}$, i.e., $\mu_S(\{x^{(i)}\})=\mu_S^{(i)} $, $ i=1,\ldots,S$.
 
Thus, by choosing $\chi_S$ and assigning a discrete probability measure $\mu_S\in\mathcal{U}_S$ to $\chi_S$, a probabilistic decision is determined, where decision variables are randomly sampled from $\chi_S$ according to $\mu_S$. For the probabilistic decision associated with $\chi_S$ and $\mu_S$, the expectations of objective and probability values are $\displaystyle\sum_{i=1}^SJ(x^{(i)})\mu_S^{(i)}$ and $\displaystyle\sum_{i=1}^S \mathbb{P}(x^{(i)})\mu_S^{(i)}$, respectively. We can optimize $\displaystyle\sum_{i=1}^SJ(x^{(i)})\mu_S^{(i)}$ under a chance constraint $\displaystyle \sum_{i=1}^S \mathbb{P}(x^{(i)})\mu_S^{(i)}\geq 1-\alpha$, where $(\mu_S,\chi_S)$ serves as decision variable. This optimization problem can be formulated as follows:
\begin{align}\label{eq:til_P_alpha_S}
    & \min_{\mu_S\in \mathcal{U}_S,\chi_S\in\mathcal{X}^S} \sum_{i=1}^SJ(x^{(i)})\mu_S^{(i)} \tag{$ \tilde{P}_{\alpha}(S)$},\\
    &{\normalfont \mathsf{s.t.}}\quad \sum_{i=1}^S \mathbb{P}(x^{(i)})\mu_S^{(i)}\geq 1-\alpha. \nonumber %
\end{align}
Define the feasible set $\tilde{\mathcal{U}}_\alpha(S)$ of Problem \ref{eq:til_P_alpha_S} by 
\begin{equation*}
    \tilde{\mathcal{U}}_\alpha(S):=\Big\{(\mu_S,\chi_S):\sum_{i=1}^S \mathbb{P}(x^{(i)})\mu_S^{(i)}\geq 1-\alpha\Big\}.
\end{equation*}
Since $\displaystyle\sum_{i=1}^S \mathbb{P}(x^{(i)})\mu_S^{(i)}: \mathcal{U}_S\times\mathcal{X}^S \to \mathbb{R} $ is a continuous function, and its domain $\mathcal{U}_S\times\mathcal{X}^S$ is compact, it follows that the feasible set $\tilde{\mathcal{U}}_\alpha(S)$ of Problem \ref{eq:til_P_alpha_S} is also compact. As a result, the optimal solution of Problem \ref{eq:til_P_alpha_S} exists. 

Define the optimal objective value $\tilde{\mathcal{J}}_\alpha(S)$ of Problem \ref{eq:til_P_alpha_S} by 
\begin{equation}
\label{eq:til_mathJ_S}
    \tilde{\mathcal{J}}_\alpha(S):=\min_{(\mu_S,\chi_S)\in  \tilde{\mathcal{U}}_\alpha(S)} \sum_{i=1}^SJ(x^{(i)})\mu_S^{(i)}.
\end{equation}
The optimal solution set of Problem \ref{eq:til_P_alpha_S} is written as 
\begin{equation*}
    \tilde{A}_{\alpha}(S):=
    \Big\{(\mu_S,\chi_S)\in  \tilde{\mathcal{U}}_\alpha(S):\sum_{i=1}^SJ(x^{(i)})\mu_S^{(i)}
    =\tilde{\mathcal{J}}_\alpha(S)\Big\}.        
\end{equation*}
For $\tilde{\mathcal{J}}_\alpha(S)$, we present the following theorem.
\begin{mytheo}\label{theo:sample_approx_convergence}
    Suppose that Assumptions \ref{assump:J_h} and \ref{assump:P_alpha_opt_sol_seq} hold. Then, we have
    \begin{equation}\label{eq:J_converge_measure}
        \lim_{S\rightarrow\infty} \tilde{\mathcal{J}}_\alpha(S)=\mathcal{J}^*_{\alpha}.
    \end{equation}
\end{mytheo}

\begin{pf}
    For $S\in\mathbb{N}$, let $\left(\mu_S,\chi_S\right)\in\tilde{\mathcal{U}}_{\alpha}(S)$ be a feasible solution to Problem \ref{eq:til_P_alpha_S} with $\chi_S=(x^{(1)},\ldots,x^{(S)})$. Let $\mu_S^{\mathsf{m}}\in M(\mathcal{X})$ be a probability measure that satisfies $\mu_S^{\mathsf{m}}\left(\{x^{(i)}\}\right)=\mu_S^{(i)}, i = 1,\ldots,S $. Then, we have
    \begin{equation*}
        \int_{\mathcal{X}}\mathbb{P}(x)\mathsf{d}\mu_S^{\mathsf{m}}(x)=\sum_{i=1}^{S}\mathbb{P}(x^{(i)})\mu_S^{(i)}\geq 1-\alpha,
    \end{equation*}
    which implies that $\mu_S^{\mathsf{m}}\in M_{\alpha}(x)$, where $M_{\alpha}(x)$ is the feasible solution set of Problem \ref{eq:P_alpha} defined by \eqref{eq:M_alpha}. Then, for all $\left(\mu_S,\chi_S\right)\in\tilde{\mathcal{U}}_{\alpha}(S),$ $S\in\mathbb{N}$, it holds that
    \begin{equation}\label{eq:J_CS_J_opt_1}
        \sum_{i=1}^S J(x^{(i)})\mu_S^{(i)}=\int_{\mathcal{X}}J(x)\mathsf{d}\mu_S^{\mathsf{m}}(x)\geq \mathcal{J}^*_{\alpha},
    \end{equation}
    which implies $\tilde{\mathcal{J}}_\alpha(S)\geq\mathcal{J}^*_{\alpha},$ for all $S\in\mathbb{N}$. Therefore, we have
    \begin{equation}\label{eq:J_CS_lim_J_opt_1}
        \liminf_{S\rightarrow\infty} \tilde{\mathcal{J}}_\alpha(S)\geq\mathcal{J}^*_{\alpha}.
    \end{equation}
    Assumption \ref{assump:P_alpha_opt_sol_seq} implies that there exists a sequence $\{\mu_k\}_{k=1}^{\infty}\subseteq M(\mathcal{X})$ that converges to an optimal solution $\mu^*$ of Problem \ref{eq:P_alpha} associated with $J(\cdot)$, such that $\int_{\mathcal{X}}\mathbb{P}(x) \mathsf{d}\mu_k> 1-\alpha,  ~~\forall  k\ge 1.$
    
    For any given $\varepsilon>0$, by the definition of convergence associated with $J(\cdot)$ (Definition \ref{def:weak_convergence}),  there exists $k_{\varepsilon}$ such that, if $k\geq k_{\varepsilon}$,   
    \begin{equation}\label{eq:J_mu_k_J_star}
       0\leq \int_{\mathcal{X}}J(x)\mathsf{d}\mu_k(x)-\mathcal{J}^*_{\alpha}\leq\frac{\varepsilon}{2},
    \end{equation}
    since $\mu^*$ is an optimal solution of Problem \ref{eq:P_alpha}, i.e.,  $\mathcal{J}^*_\alpha=\int_{\mathcal{X}}J(x)\mathsf{d}\mu^*$ by \eqref{eq:mathJ_star}.

    Let $\chi_{S}^{k_{\varepsilon}}:= (x^{(1)}_{k_{\varepsilon}},\ldots,x^{(S)}_{k_{\varepsilon}}), S\in \mathbb{N}$ be a sample set obtained by sampling from $\mathcal{X}$ according to probability measure $\mu_{k_{\varepsilon}}$. By law of large numbers \cite[p. 457]{Shapiro}, for any continuous function $f\in \mathscr{C}(\mathcal{X})$, as $S\rightarrow\infty$, w. p. 1, we have 
    \begin{equation}
        \label{eq:LLN}
        \frac{1}{S}\sum_{i=1}^{S} f(x^{(i)}_{k_{\varepsilon}})\rightarrow \mathbb{E}_{x\sim\mu_{k_{\varepsilon}}}\left\{f(x)\right\}=\int_{\mathcal{X}}f(x)\mathsf{d}\mu_{k_{\varepsilon}}.
    \end{equation}
    Notice that $J(\cdot)$ and $\mathbb{P}(\cdot)$ are continuous functions. Then, \eqref{eq:LLN} also holds by replacing $f(\cdot)$ by either $J(\cdot)$ or $\mathbb{P}(\cdot)$. As a result, for given 
    \begin{equation}
        \label{eq:vep1}
        \varepsilon_1=\min\left\{\frac{\varepsilon}{2}, \int_{\mathcal{X}}\mathbb{P}(x) \mathsf{d}\mu_{k_{\varepsilon}}-( 1-\alpha) \right\}>0,
    \end{equation}
    there exists $S(\varepsilon)=S(\varepsilon, k_{\varepsilon})$ such that, if $S\geq S(\varepsilon)$, w. p. 1, the following conditions hold:
    \begin{align}
        \label{eq:P_int_approx}
       & \left|\frac{1}{S}\sum_{i=1}^{S} \mathbb{P}({x}^{(i)}_{{k_{\varepsilon}}})-\int_{\mathcal{X}}\mathbb{P}(x)\mathsf{d}\mu_{k_{\varepsilon}}\right|\leq \varepsilon_1,\\
        \label{eq:J_int_approx}
        &\left|\frac{1}{S}\sum_{i=1}^{S} J(\tilde{x}^{(i)}_{k_{\varepsilon}})-\int_{\mathcal{X}}J(x)\mathsf{d}\mu_{k_{\varepsilon}}\right|\leq \varepsilon_1,
    \end{align}
    for any $\chi_{S}^{k_{\varepsilon}}.$ 
    Define a  discrete probability measure $\mu^{S}_k\in M(\mathcal{X}), S\geq S(\varepsilon)$ such that $\mu_S^{k_{\varepsilon}(i)} =\frac{1}{S},\ i = 1,\ldots, S,$
    where $\mu_S^{k_{\varepsilon}(i)}$ is the discrete probability assigned to ${x}^{(i)}_{k_{\varepsilon}}$, i.e., 
    $\mu_S^{k_{\varepsilon}}(\{{x}^{(i)}_{k_{\varepsilon}}\})=
    \mu_S^{k_{\varepsilon}(i)}$. 
    From \eqref{eq:vep1}, and \eqref{eq:P_int_approx}, we have
    \begin{equation*}
       \sum_{i=1}^{S}   \mathbb{P}({x}^{(i)}_{k_{\varepsilon}}) \mu_S^{k_{\varepsilon} (i)}\geq \int_{\mathcal{X}}\mathbb{P}(x)\mathsf{d}\mu_{k_{\varepsilon}}-{\varepsilon}_1\geq 1-\alpha,
    \end{equation*}
    which implies for each  $S\geq S(\varepsilon)$,  w. p. 1,
    $(\mu_S^{k_{\varepsilon}}, \chi^{k_{\varepsilon}}_{S})$ is a feasible solution of Problem \ref{eq:til_P_alpha_S} and thus, w. p. 1, 
    \begin{equation}
        \label{eq:mu_S_k_cost_com}
        \sum_{i=1}^{S}   J({x}_{k_{\varepsilon}}^{(i)}) \mu_S^{k_{\varepsilon}(i)} \geq\tilde{\mathcal{J}}_\alpha(S).
    \end{equation}
    Considering \eqref{eq:J_int_approx} and \eqref{eq:mu_S_k_cost_com}, we have\footnote{From \eqref{eq:J_int_approx} and \eqref{eq:mu_S_k_cost_com}, we firstly obtain the inequality 
    $$\tilde{\mathcal{J}}_\alpha(S)-\int_{\mathcal{X}}J(x)\mathsf{d}\mu^{k_{\varepsilon}}\leq\tilde{\varepsilon}_1,  \text{w.p.1.}$$
    hold. But noticing  the both of $\tilde{\mathcal{J}}_\alpha(S)$, and $\int_{\mathcal{X}}J(x)\mathsf{d}\mu^{k_{\varepsilon}}$ are deterministic value, we can remove ``w. p. 1.'' above inequality, and obtain \eqref{eq:J_CS_J_mu_k}.
    }
    \begin{equation}
        \label{eq:J_CS_J_mu_k}
        \tilde{\mathcal{J}}_\alpha(S)-\int_{\mathcal{X}}J(x)\mathsf{d}\mu^{k_{\varepsilon}}\leq\tilde{\varepsilon}_1.
    \end{equation}
    For $S\geq S(\varepsilon)$, combining \eqref{eq:J_mu_k_J_star} and \eqref{eq:J_CS_J_mu_k},  we have
    \begin{equation}
        \label{eq:J_CS_J_star}
        \tilde{\mathcal{J}}_\alpha(S)-\mathcal{J}^*_\alpha\leq\frac{\varepsilon}{2}+\varepsilon_1\le \varepsilon,
    \end{equation}
    which implies $\displaystyle \limsup_{S\rightarrow\infty} \tilde{\mathcal{J}}_\alpha(S)\leq\mathcal{J}^*_{\alpha}+\varepsilon$. 
    
    Finally, by arbitrary  of $\varepsilon$, we have 
    \begin{equation}
        \label{eq:J_CS_lim_J_opt_2}
        \limsup_{S\rightarrow\infty} \tilde{\mathcal{J}}_\alpha(S)\leq\mathcal{J}^*_{\alpha}.
    \end{equation}
   \eqref{eq:J_CS_lim_J_opt_1} and \eqref{eq:J_CS_lim_J_opt_2} forms \eqref{eq:J_converge_measure} holds. \qed
\end{pf}

From Theorem \ref{theo:sample_approx_convergence}, it implies that if the $S\rightarrow\infty$, the optimal solution $\tilde{\mu}_S\in\tilde{A}_{\alpha}(S)$ of Problem \ref{eq:til_P_alpha_S} can be used as an approximate optimal solution of Problem \ref{eq:P_alpha}. 
However, solving Problem \ref{eq:til_P_alpha_S} becomes computationally intractable when $S\rightarrow\infty$. 
Theorem \ref{theo:P_alpha_2_point_approx} show that the optimal objective value $\tilde{\mathcal{J}}_\alpha(S)$ equals $\mathcal{J}^*_\alpha$ when $S=2$. 
We further give some preparations for Proving Theomre \ref{theo:P_alpha_2_point_approx}.


For a given integer $S\in\mathbb{N}$, define a multi-sample realization of violation probabilities as
\begin{equation}\label{eq:vio_pro_set_prop_1}
    \mathcal{E}_S=(\tilde{\alpha}^{(1)},\ldots,\tilde{\alpha}^{(S)}),
\end{equation}
where $\tilde{\alpha}^{(i)}=1-\mathbb{P}(x^{(i)})$. 
The multi-sample realization $\mathcal{E}_S$ is an element in the space $[0,1]^S.$
Note that each $\tilde{\alpha}^{(i)}$ is a threshold of violation probability in Problem $Q_{\tilde{\alpha}^{(i)}}$ (Problem \ref{eq:Q_alpha} when $\alpha=\tilde{\alpha}^{(i)}$). 
For a given $\mathcal{E}_S$, we can construct a corresponding optimal objective value set $\{J^*_{\tilde{\alpha}^{(1)}},\ldots,J^*_{\tilde{\alpha}^{(S)}}\}$,
where $J^*_{\tilde{\alpha}^{(i)}}$ is the optimal objective value of Problem $Q_{\tilde{\alpha}^{(i)}}$. 
Let $\mathcal{V}_{S}:=\big\{\nu_{S}\in[0,1]^S:\displaystyle\sum_{i=1}^S\nu_{S}^{(i)}=1\big\}$ be a set of discrete probability measures that defined on $\mathcal{E}_S$. 
Define a decision variable $\bm{\theta}_S:=(\nu_S,\mathcal{E}_S)\in\Theta_{S}:=\mathcal{V}_{S}\times [0,1]^S.$ 
An optimization problem of selecting $\bm{\theta}_S$ is formulated by 
\begin{align} 
    &\min_{\bm{\theta}_S\in\Theta_{S}} \,\, \sum_{i=1}^SJ^*_{\tilde{\alpha}^{(i)}}\nu_S^{(i)}, \tag{$\tilde{V}_{\alpha}(S)$} \label{eq:til_V_alpha_S}\\
    &\mathsf{s.t.}\quad \sum_{i=1}^{S}(1-\tilde{\alpha}^{(i)})\nu_S^{(i)}\geq 1-\alpha. \nonumber
\end{align} 
Define the feasible set $\Theta_{S,\alpha}$ of \ref{eq:til_V_alpha_S} by
\begin{equation*}
    \Theta_{S,\alpha}:=
    \Big\{\bm{\theta}_{S}\in\Theta_{S}:\sum_{i=1}^{S}(1-\tilde{\alpha}^{(i)})\nu_{S}^{(i)}\geq 1-\alpha\Big\}.
\end{equation*}
Unlike Problem \ref{eq:til_P_alpha_S}, Problem \ref{eq:til_V_alpha_S} optimizes finite combinations of violation probability thresholds to achieve the minimal mean objective value while satisfying chance constraints. The relationship between Problems \ref{eq:til_P_alpha_S} and \ref{eq:til_V_alpha_S} is established in Proposition \ref{prop:A_nu_alpha}.

\begin{myprop}
    \label{prop:A_nu_alpha}
    Suppose Assumption \ref{assump:J_h} holds. For all $S\in\mathbb{N}$, there exists a feasible solution $\bm{\theta}_{S}:=(\nu_{S}, \tilde{\alpha}^{(1)},\ldots,\tilde{\alpha}^{(S)})\in\Theta_{S,\alpha}$ of Problem \ref{eq:til_V_alpha_S} that satisfies 
    \begin{equation}
        \label{eq:mu_nu_alpha}
        \sum_{i=1}^SJ^*_{\tilde{\alpha}^{(i)}}\nu_S^{(i)}=\tilde{\mathcal{J}}_\alpha(S),
    \end{equation}
    where $\tilde{\mathcal{J}}_\alpha(S)$ is the optimal objective value of Problem \ref{eq:til_P_alpha_S} defined by \eqref{eq:til_mathJ_S}.
\end{myprop}

\begin{pf}
    For a given integer $S,$ define an optimal solution of Problem \ref{eq:til_P_alpha_S} by $\left(\tilde{\mu}_{\alpha},\tilde{\chi}_S\right)\in\tilde{A}_{\alpha}(S).$ 
    Here, we have $\tilde{\chi}_S=(x^{(1)},\ldots,x^{(S)}).$ 
    Since $\left(\tilde{\mu}_{\alpha},\tilde{\chi}_S\right)$ is a feasible solution, we have 
    \begin{equation*}
        \sum_{i=1}^S J(x^{(i)})\tilde{\mu}_{\alpha}^{(i)}=\tilde{\mathcal{J}}_\alpha(S),\ \sum_{i=1}^S \mathbb{P}(x^{(i)})\tilde{\mu}_{\alpha}^{(i)}\geq 1-\alpha.
    \end{equation*}
    Choose $\bm{\theta}_S=(\nu_S,\mathcal{E}_S)$ satisfying $\nu_S^{(i)}=\tilde{\mu}_\alpha^{(i)}, i = 1,\ldots,S$ and $\tilde{\alpha}^{(i)}=1-\mathbb{P}(x^{(i)}),$ which further implies  
    \begin{equation*}
        \sum_{i=1}^S (1-\tilde{\alpha}^{(i)})\nu_S^{(i)}=\sum_{i=1}^S \mathbb{P}(x^{(i)})\tilde{\mu}_\alpha^{(i)}\geq 1-\alpha.
    \end{equation*}
    Therefore, $\bm{\theta}_S=(\nu_S,\mathcal{E}_S)$ is a feasible solution of Problem \ref{eq:til_V_alpha_S} and then it holds that 
    \begin{align}
        \label{eq:calJ_E_plus_S}
        \sum_{i=1}^S J^*_{\tilde{\alpha}^{(i)}}\nu_S^{(i)}\leq \sum_{i=1}^S J(x^{(i)})\tilde{\mu}_\alpha^{(i)}=\tilde{\mathcal{J}}_\alpha(S).
    \end{align} 
    To show $\displaystyle\sum_{i=1}^S J^*_{\tilde{\alpha}^{(i)}}\nu_S^{(i)}=\tilde{\mathcal{J}}_\alpha(S)$, we need to prove that $\displaystyle\sum_{i=1}^S J^*_{\tilde{\alpha}^{(i)}}\nu_S^{(i)}\geq\tilde{\mathcal{J}}_\alpha(S)$ also holds.
    
    For $\mathcal{E}_S$ defined in \eqref{eq:vio_pro_set_prop_1}, define a multiple-sample realization of decision variables as $\hat{\chi}_S=(\hat{x}^{(1)}, \ldots, \hat{x}^{(S)}).$ 
    Here, $\hat{x}^{(i)}\in X_{\tilde{\alpha}^{(i)}}$ is one of the optimal solution of Problem $Q_{\tilde{\alpha}^{(i)}}.$ 
    The optimality implies $J(\hat{x}^{(i)})=J^*_{\tilde{\alpha}^{(i)}}.$
    Besides, $\hat{x}^{(i)}$ is a feasible solution of Problem $Q_{\tilde{\alpha}^{(i)}},$ and we thus have $\mathbb{P}(\hat{x}^{(i)})\geq 1-\tilde{\alpha}^{(i)}.$ 
    
    Choosing $\hat{\mu}_S^{(i)}=\nu_S^{(i)}$ and $\left(\hat{\mu}_S,\hat{\chi}_S\right)$ is then a feasible solution of Problem \ref{eq:til_P_alpha_S} since we have
    \begin{equation*}
        \sum_{i=1}^S \mathbb{P}(\hat{x}^{(i)})\hat{\mu}_{S}^{(i)}\geq \sum_{i=1}^S (1-\tilde{\alpha}^{(i)})\nu_S^{(i)} \geq 1-\alpha.
    \end{equation*} 
    Furthermore, due to the feasibility of $\left(\hat{\mu}_S,\hat{\chi}_S\right)$, we have that its objective value is larger than the optimal one, which can be written by
    \begin{align}
        \label{eq:calJ_C_plus_S}
        \sum_{i=1}^S J^*_{\tilde{\alpha}^{(i)}}\nu_{S}^{(i)}=\sum_{i=1}^S J(\hat{x}^{(i)})\hat{\mu}_{S}^{(i)}\geq\tilde{\mathcal{J}}_\alpha(S).
    \end{align}
    From \eqref{eq:calJ_E_plus_S} and \eqref{eq:calJ_C_plus_S}, we can conclude \eqref{eq:mu_nu_alpha}. \qed 
\end{pf}

Proposition \ref{prop:A_nu_alpha} can be used to obtain the following lemma. \begin{mylemma} \label{lemma:V_alpha_2_point_approx} Suppose Assumptions \ref{assump:J_h}-\ref{assump:P_alpha_opt_sol_seq} both hold. There exists a feasible solution $\bm{\theta}_{\mathsf{m}}:=(\nu_{\mathsf{m}}^{(1)}, \nu_{\mathsf{m}}^{(2)}, \tilde{\alpha}_{\mathsf{m}}^{(1)}, \tilde{\alpha}_{\mathsf{m}}^{(2)})\in\Theta_{2,\alpha}$ of Problem $\tilde{V}_\alpha(2)$ that satisfies \begin{equation} \label{eq:mu_nu_alpha_opt} \sum_{i=1}^2J^*_{\tilde{\alpha}_{\mathsf{m}}^{(i)}}\nu_{\mathsf{m}}^{(i)}=\mathcal{J}^*_\alpha. \end{equation} \end{mylemma}
\begin{pf}
    Denote 
    $\mathcal{H} := \left\{(\tilde{\alpha}, J^*_{\tilde{\alpha}}):\tilde{\alpha}\in(0,1]\right\}$,
    where $\tilde{\alpha}$ represents a violation probability threshold and $J^*_{\tilde{\alpha}}$ is the optimal objective value of Problem $Q_{\tilde{\alpha}}$. 
    Besides, use $\mathsf{conv}(\cdot)$ to represent the convex hull of a set. 
    Consider the following optimization problem:
    \begin{align} 
        &\underset{(\tilde{\alpha}_{\mathsf{h},\alpha},J_{\mathsf{h}})\in \mathsf{conv}(\mathcal{H})} {{\normalfont\mathsf{min}}} \,\, J_{\mathsf{h}}, \tag{$H_{\alpha}$}\\
        &{\normalfont \mathsf{s.t.}}\quad 1-\tilde{\alpha}_{\mathsf{h},\alpha}\geq 1-\alpha. \nonumber
    \end{align}
    Under Assumptions \ref{assump:J_h}-\ref{assump:P_alpha_opt_sol_seq}, the set $\mathcal{H}$ is bounded in $\mathbb{R}^2$ and the mapping $\tilde{\alpha}\mapsto J^*_{\tilde{\alpha}}$ is well-defined on $(0,1]$. Hence, $\mathsf{conv}(\mathcal{H})$ is a bounded convex set. Moreover, by Proposition \ref{prop:A_nu_alpha}, for any $S\in\mathbb{N}$ there exists a feasible point $\bm{\theta}_S$ of Problem \ref{eq:til_V_alpha_S}, and for such a feasible point we have $\sum_{i=1}^S\tilde{\alpha}^{(i)}\nu_S^{(i)}\leq \alpha$. Therefore, the point $\left(\sum_{i=1}^S\tilde{\alpha}^{(i)}\nu_S^{(i)},\sum_{i=1}^S J^*_{\tilde{\alpha}^{(i)}}\nu_S^{(i)}\right)$ belongs to $\mathsf{conv}(\mathcal{H})$ and satisfies the constraint of Problem $H_\alpha$, which implies the feasible set of Problem $H_\alpha$ is nonempty. Since the objective is linear in $J_{\mathsf{h}}$ and the feasible set is bounded, Problem $H_\alpha$ admits an optimal solution. Let $(\tilde{\alpha}_{\mathsf{h},\alpha}^{\diamond},J^\diamond_{\mathsf{h},\alpha})$ be an optimal solution of Problem $H_\alpha$. 
    
    \textbf{Step 1: Establishing $J^\diamond_{\mathsf{h},\alpha} \leq \mathcal{J}^*_\alpha$}.  
    For any $S\in\mathbb{N}$, let $\bm{\theta}_S := \{\nu_S, \tilde{\alpha}^{(1)},\ldots,\tilde{\alpha}^{(S)}\}$ be a feasible solution of Problem \ref{eq:til_V_alpha_S} satisfying \eqref{eq:mu_nu_alpha} in Proposition \ref{prop:A_nu_alpha}. Define  
    \begin{equation*}
        \tilde{\alpha}_{\mathsf{mean}}(\bm{\theta}_S):=\sum_{i=1}^S\tilde{\alpha}^{(i)}\nu_S^{(i)},\  J^*_{\mathsf{mean}}(\bm{\theta}_S):=\sum_{i=1}^S J^*_{\tilde{\alpha}^{(i)}}\nu_S^{(i)}.
    \end{equation*}
    
    By the definition of convex hull, $\left(\tilde{\alpha}_{\mathsf{mean}}(\bm{\theta}_S),J^*_{\mathsf{mean}}(\bm{\theta}_S)\right) \in \mathsf{conv}(\mathcal{H})$. Since $\bm{\theta}_S$ is feasible, we have $1-\tilde{\alpha}_{\mathsf{mean}}(\bm{\theta}_S) \geq 1-\alpha$ and $J^*_{\mathsf{mean}}(\bm{\theta}_S) = \tilde{\mathcal{J}}_\alpha(S)$. Thus, for all $S \in \mathbb{N}$,
    \begin{equation*}
        J^\diamond_{\mathsf{h},\alpha} \leq \tilde{\mathcal{J}}_\alpha(S).
    \end{equation*}
    Taking $S \to \infty$ and applying Theorem \ref{theo:sample_approx_convergence}, we obtain
    \begin{equation*}
        J^\diamond_{\mathsf{h},\alpha} \leq \mathcal{J}^*_\alpha.
    \end{equation*}
    
    \textbf{Step 2: Establishing the Boundary Condition}.  
    We show that $(\tilde{\alpha}_{\mathsf{h},\alpha}^{\diamond},J^\diamond_{\mathsf{h},\alpha})$ cannot be an interior point of $\mathsf{conv}(\mathcal{H})$. Suppose for contradiction that $(\tilde{\alpha}_{\mathsf{h},\alpha}^{\diamond},J^\diamond_{\mathsf{h},\alpha})$ is an interior point of $\mathsf{conv}(\mathcal{H})$. Then, there exists $\varepsilon > 0$ such that $\mathcal{B}_{\varepsilon}((\tilde{\alpha}_{\mathsf{h},\alpha}^{\diamond}, J^\diamond_{\mathsf{h},\alpha})) \subset \mathsf{conv}(\mathcal{H})$. For any $\tilde{\varepsilon} < \varepsilon$, the point $(\tilde{\alpha}_{\mathsf{h},\alpha}^{\diamond}, J^\diamond_{\mathsf{h},\alpha} - \tilde{\varepsilon})$ also belongs to $\mathsf{conv}(\mathcal{H})$ and still satisfies $1-\tilde{\alpha}_{\mathsf{h},\alpha}^{\diamond} \geq 1-\alpha$. However, this contradicts the optimality of $J^\diamond_{\mathsf{h},\alpha}$, implying that $(\tilde{\alpha}_{\mathsf{h},\alpha}^{\diamond},J^\diamond_{\mathsf{h},\alpha})$ must lie on the boundary of $\mathsf{conv}(\mathcal{H})$.
    
    \textbf{Step 3: Application of Carath\'{e}odory's Theorem}.  
    Since $(\tilde{\alpha}_{\mathsf{h},\alpha}^{\diamond},J^\diamond_{\mathsf{h},\alpha})$ is a boundary point of the convex set $\mathsf{conv}(\mathcal{H})$, by the supporting hyperplane theorem \cite[p. 133]{Luenberger}, there exists a line $\mathcal{L}$ passing through $(\tilde{\alpha}_{\mathsf{h},\alpha}^{\diamond},J^\diamond_{\mathsf{h},\alpha})$ that contains $\mathsf{conv}(\mathcal{H})$ in one of its closed half-spaces. Since $\mathcal{L}$ is a one-dimensional linear space, we conclude that $(\tilde{\alpha}_{\mathsf{h},\alpha}^{\diamond}, J^\diamond_{\mathsf{h},\alpha}) \in \mathsf{conv}(\mathcal{L} \cap \mathcal{H}).$
    By Carath\'{e}odory's theorem \cite{Eckhoff}, this implies that $(\tilde{\alpha}_{\mathsf{h},\alpha}^{\diamond},J^\diamond_{\mathsf{h},\alpha})$ is a convex combination of at most two points in $\mathcal{L} \cap \mathcal{H}$. That is, there exist weights $\nu_{\mathsf{m}}^\diamond \in \mathcal{V}_2$ and two values $(\tilde{\alpha}^{(1)}_{\mathsf{m}}, \tilde{\alpha}^{(2)}_{\mathsf{m}}) \in [0,1]^2$ such that  
    $J^\diamond_{\mathsf{h},\alpha} = J^*_{\tilde{\alpha}^{(1)}_{\mathsf{m}}} \nu_{\mathsf{m}}^{\diamond(1)} + J^*_{\tilde{\alpha}^{(2)}_{\mathsf{m}}} \nu_{\mathsf{m}}^{\diamond(2)},
        \tilde{\alpha}_{\mathsf{h},\alpha}^{\diamond} = \tilde{\alpha}^{(1)}_{\mathsf{m}} \nu_{\mathsf{m}}^{\diamond(1)} + \tilde{\alpha}^{(2)}_{\mathsf{m}} \nu_{\mathsf{m}}^{\diamond(2)}.$
    
    Since $1 - \tilde{\alpha}_{\mathsf{h},\alpha}^{\diamond} \geq 1-\alpha$, it follows that $(\nu_{\mathsf{m}}^{\diamond(1)},\nu_{\mathsf{m}}^{\diamond(2)},\tilde{\alpha}_\mathsf{m}^{(1)},\tilde{\alpha}_\mathsf{m}^{(2)})$ is a feasible solution of Problem \ref{eq:til_V_alpha_S} for $S=2$. Defining optimal solutions $\tilde{x}_{\mathsf{m}}^{(1)}\in X_{\tilde{\alpha}_\mathsf{m}^{(1)}}$ and $\tilde{x}_{\mathsf{m}}^{(2)}\in X_{\tilde{\alpha}_\mathsf{m}^{(2)}}$ for $Q_{\tilde{\alpha}_\mathsf{m}^{(1)}}$ and $Q_{\tilde{\alpha}_\mathsf{m}^{(2)}}$, we obtain
    \begin{align*}
        J^\diamond_{\mathsf{h},\alpha} = \sum_{i=1}^2 J(\tilde{x}_\mathsf{m}^{(i)}) \nu_{\mathsf{m}}^{\diamond(i)}.
    \end{align*}
    By feasibility of $(\nu_{\mathsf{m}}^{\diamond(1)},\nu_{\mathsf{m}}^{\diamond(2)},\tilde{\alpha}_\mathsf{m}^{(1)},\tilde{\alpha}_\mathsf{m}^{(2)})$ for Problem \ref{eq:til_V_alpha_S} with $S=2$, and by the definition of $\tilde{\mathcal{J}}_\alpha(2)$ as the optimal objective value of Problem $\tilde{V}_\alpha(2)$, we have
    \begin{align*}
        J^\diamond_{\mathsf{h},\alpha} \geq \tilde{\mathcal{J}}_\alpha(2).
    \end{align*}
    Moreover, since $\tilde{V}_\alpha(2)$ is a restriction of the corresponding measure formulation, its optimal value satisfies $\tilde{\mathcal{J}}_\alpha(2) \geq \mathcal{J}^*_\alpha$. Therefore,
    \begin{align*}
        J^\diamond_{\mathsf{h},\alpha} \geq \tilde{\mathcal{J}}_\alpha(2) \geq \mathcal{J}^*_\alpha.
    \end{align*}
    Combining this with \textbf{Step 1}, we obtain $J^\diamond_{\mathsf{h},\alpha} = \mathcal{J}^*_\alpha$. Hence, there exists a feasible solution $\bm{\theta}_{\mathsf{m}}:=(\nu_{\mathsf{m}}^{(1)}, \nu_{\mathsf{m}}^{(2)}, \tilde{\alpha}_{\mathsf{m}}^{(1)}, \tilde{\alpha}_{\mathsf{m}}^{(2)})\in\Theta_{2,\alpha}$ of Problem $\tilde{V}_\alpha(2)$ that satisfies \eqref{eq:mu_nu_alpha_opt}. \qed
\end{pf}

\section{Proof of Theorem \ref{theo:P_alpha_2_point_approx}} \label{proof:theo_P_alpha_2_point_approx}

Based on Theorem \ref{theo:sample_approx_convergence} and Lemma \ref{lemma:V_alpha_2_point_approx}, the proof of Theorem~\ref{theo:P_alpha_2_point_approx} is summarized as follows:

By \eqref{eq:J_CS_J_opt_1} in the Proof of Theorem \ref{theo:sample_approx_convergence}, we have $\mathcal{J}^*_\alpha\leq\mathcal{J}^{\mathsf{w}}_{\alpha}$, where $\mathcal{J}^{\mathsf{w}}_{\alpha}=\tilde{J}_\alpha(S)$ with $S=2$. Then, we need to show that $\mathcal{J}^*_\alpha\geq\mathcal{J}^{\mathsf{w}}_{\alpha}$. By Lemma \ref{lemma:V_alpha_2_point_approx}, we can find a feasible solution $\bm{\theta}_\mathsf{m}:=\{\nu_\mathsf{m}^{(1)},\nu_\mathsf{m}^{(2)},\tilde{\alpha}_\mathsf{m}^{(1)},\tilde{\alpha}_\mathsf{m}^{(2)}\}\in\Theta_{2,\alpha}$ of Problem \ref{eq:til_V_alpha_S} with $S=2$ that satisfies \eqref{eq:mu_nu_alpha_opt}, namely, $\mathcal{J}^*_\alpha=\displaystyle\sum_{i=1}^2 J_{\tilde{\alpha}^{(i)}_{\mathsf{m}}}^*\nu_{\mathsf{m}}^{(i)}$. Choose $\mu_{\mathsf{m}}=\nu_{\mathsf{m}}$ and $x^{*(1)} \in X_{\tilde{\alpha}^{(1)}_{\mathsf{m}}},\ x^{*(2)}\in X_{\tilde{\alpha}^{(2)}_{\mathsf{m}}}$. Since $\bm{\theta}_\mathsf{m}$ is a feasible solution of Problem $\tilde{V}_\alpha(2)$, we have
\begin{equation}\label{eq:tilde_mu_opt_m_P}
    \sum_{i=1}^2 \mathbb{P}(x^{*(i)})\mu_{\mathsf{m}}^{(i)}\geq\sum_{i=1}^2 \left(1-\tilde{\alpha}^{(i)}_{\mathsf{m}}\right)\nu_{\mathsf{m}}^{(i)}\geq 1-\alpha.
\end{equation}
Notice that \eqref{eq:tilde_mu_opt_m_P} implies that $(\mu_{\mathsf{m}}^{(1)},\mu_{\mathsf{m}}^{(2)},x^{*(1)},x^{*(2)})$ is a feasible solution of Problem \ref{eq:W_alpha}. By Lemma \ref{lemma:V_alpha_2_point_approx}, we have
\begin{equation*}
    \mathcal{J}^*_\alpha=\sum_{i=1}^2 J_{\tilde{\alpha}^{(i)}_{\mathsf{m}}}^*\nu_{\mathsf{m}}^{(i)}=\sum_{i=1}^2 J(x^{*(i)})\mu_{\mathsf{m}}^{(i)} \geq \mathcal{J}^{\mathsf{w}}_\alpha.
\end{equation*}

Thus, we conclude $\mathcal{J}^{\mathsf{w}}_{\alpha}=\mathcal{J}^*_\alpha$.\qed

\section{Proof of Corollary \ref{coro:z_m_opt}}
\label{proof:coro_z_m_opt}

Suppose that $x^{*(1)}\notin X_{\tilde{\alpha}^{(1)}}$. Then, $J(x^{*(1)})>J^*_{\tilde{\alpha}^{(1)}}$ holds. Let $\hat{x}^{(1)}\in X_{\tilde{\alpha}^{(1)}}$. Then, $\mathbb{P}(\hat{x}^{(1)})\geq 1-\tilde{\alpha}^{(1)}=\mathbb{P}(x^{*(1)})$ and $J(\hat{x}^{(1)})=J^*_{\tilde{\alpha}^{(1)}}<J(x^{*(1)})$. We have
 \begin{align}
    &\displaystyle \mathbb{P}(\hat{x}^{(1)})\mu_{\mathsf{m}}^{*(1)}+\mathbb{P}(x^{*(2)})\mu_{\mathsf{m}}^{*(2)}\geq\sum_{i=1}^2\mathbb{P}(x^{*(i)})\mu_{\mathsf{m}}^{*(i)} \geq 1-\alpha,\nonumber \\
    &\displaystyle J(\hat{x}^{(1)})\mu_{\mathsf{m}}^{*(1)}+J(x^{*(2)})\mu_{\mathsf{m}}^{*(2)}<\sum_{i=1}^2 J(x^{*(i)})\mu_{\mathsf{m}}^{*(i)}.\nonumber
\end{align}
Thus, $\hat{\bm{z}}_{\mathsf{m}}=(\mu_{\mathsf{m}}^{*(1)},\mu_{\mathsf{m}}^{*(2)},\hat{x}^{(1)},x^{*(2)})$ is a feasible solution of Problem \ref{eq:W_alpha} with a smaller objective function value than $\bm{z}_{\mathsf{m}}^*$, which contradicts to the optimality of $\bm{z}_{\mathsf{m}}^*\in D_\alpha$. Consequently, we conclude that $x^{*(1)}\in X_{\tilde{\alpha}^{(1)}}$. Following the above procedures, we can similarly prove that $x^{*(2)}\in X_{\tilde{\alpha}^{(2)}}$. \qed

\section{Proof of Corollary \ref{coro:W_alpha_Q_alpha}}\label{proof:coro_W_alpha_Q_alpha}

Let $\bm{z}_{\mathsf{m}}^*=(\mu_{\mathsf{m}}^{*(1)},\mu_{\mathsf{m}}^{*(2)},x^{*(1)},x^{*(2)} )\in D_\alpha$ be an optimal solution of Problem \ref{eq:W_alpha}. Let $\tilde{\alpha}^{(1)}=1-\mathbb{P}(x^{*(1)})$ and $\tilde{\alpha}^{(2)}=1-\mathbb{P}(x^{*(2)})$.

By Corollary \ref{coro:z_m_opt}, we have $J(x^{*(1)})=J^*_{\tilde{\alpha}^{(1)}}$ and $J(x^{*(2)})=J^*_{\tilde{\alpha}^{(2)}}$. Without loss of generality, suppose $\tilde{\alpha}^{(1)}<\tilde{\alpha}^{(2)}$. Let $\tilde{\alpha}^\mathsf{com}:=\mu_{\mathsf{m}}^{*(1)}\tilde{\alpha}^{(1)}+\mu_{\mathsf{m}}^{*(2)}\tilde{\alpha}^{(2)}$. Note that $\tilde{\alpha}^\mathsf{com}\leq\alpha$ since $\bm{z}_{\mathsf{m}}^*$ is feasible for \ref{eq:W_alpha}. Because $J^*_{\tilde{\alpha}}$ is a convex function of $\tilde{\alpha}$, we have
\begin{equation*}
    J^*_{\tilde{\alpha}^\mathsf{com}}\leq J^*_{\tilde{\alpha}^{(1)}}\mu_{\mathsf{m}}^{*(1)}+J^*_{\tilde{\alpha}^{(2)}}\mu_{\mathsf{m}}^{*(2)}=\mathcal{J}^*_{\alpha}.
\end{equation*}
The last equality holds due to Theorem \ref{theo:P_alpha_2_point_approx}. Since $J^*_{\tilde{\alpha}}$ is monotonically decreasing function of $\tilde{\alpha}$, we have
\begin{equation*}
    J^*_\alpha\leq J^*_{\tilde{\alpha}^\mathsf{com}}\leq\mathcal{J}^*_{\alpha}.
\end{equation*}
By \eqref{eq:J_m_mu} and \eqref{eq:P_m_mu}, $J^*_\alpha\geq \mathcal{J}^*_{\alpha}$ holds. Thus, $\mathcal{J}^*_{\alpha}=J^*_\alpha$ holds.

If $\alpha=0$, $\tilde{\alpha}^{(1)}=\tilde{\alpha}^{(2)}$ holds when $\mu_{\mathsf{m}}^{*(1)}$ and $\mu_{\mathsf{m}}^{*(2)}$ are positive, then $J(x^{*(1)})=J(x^{*(2)})=J^{*}_0$ since $x^{*(1)},x^{*(2)}\in X_0$. Thus, $\mathcal{J}_0^*=(\mu_{\mathsf{m}}^{*(1)}+\mu_{\mathsf{m}}^{*(2)})J^*_{0}=J^*_{0}$. \qed

\section{Proof of Theorem \ref{theo:uniform_convergence_W}}
\label{proof:theo_uniform_convergence_W}

By Assumption \ref{assump:opti_solu_Q_alpha}, the set $X_{\tilde{\alpha}}$ is nonempty for any $\tilde{\alpha}\in[0,1]$ and there exists $x\in\mathcal{X}$ such that $\mathbb{P}(x)> 1-\tilde{\alpha}$. By \cite[Theorem 3.4]{Pena:2020}, $\tilde{\mathbb{P}}_N^{\gamma,\varepsilon}(x)$ converges to $\mathbb{P}(x)$ as $N\rightarrow\infty$, thus there exists a sufficiently small $\varepsilon_0$ and a sufficiently large $N_0$ such that $\tilde{\mathbb{P}}_N^{\gamma,\varepsilon}(x)\geq 1-\tilde{\alpha},$ for any $ N>N_0$ w. p. 1. Notice that $\tilde{\mathbb{P}}_N^{\gamma,\varepsilon}(x)$ is continuous in $x$ and $\mathcal{X}$ is compact, and then the feasible set of Problem $\tilde{W}_{\tilde{\alpha}}^{\gamma,\varepsilon}(\mathcal{D}_N)$ is compact. Thus, $\tilde{X}_{\tilde{\alpha}}^{\gamma,\varepsilon}$ is nonempty for all $N\geq N_0$ and $\varepsilon\leq\varepsilon_0$. 

Define two sequences $\{N_k\}_{k=1}^\infty\geq N_0$ and $\{\varepsilon_k\}_{k=1}^{\infty}$, where $N_k\rightarrow\infty$ and $\varepsilon_k\rightarrow 0$. For each $N_k$ and $\varepsilon_k$, there exists a corresponding optimal solution $\tilde{\bm{z}}_{\mathsf{m},\alpha}^{\gamma,\varepsilon_k}(\mathcal{D}_{N_k})\in\tilde{D}_{\alpha}^{\gamma,\varepsilon_k}(\mathcal{D}_{N_k})$ of Problem $\tilde{W}^{\gamma,\varepsilon_k}_{\alpha}(\mathcal{D}_{N_k})$, which is represented by $\tilde{\bm{z}}_{\mathsf{m},\alpha}^k$ in this proof for simplicity. Then, we have a sequence $\{\tilde{\bm{z}}^k_{\mathsf{m},\alpha}\}_{k=1}^{\infty}$. Specify the components of $\tilde{\bm{z}}_{\mathsf{m},\alpha}^k$ by $\tilde{\bm{z}}^k_{\mathsf{m},\alpha}=(\tilde{\mu}_{\mathsf{m}}^{k(1)},\tilde{\mu}_{\mathsf{m}}^{k(2)},\tilde{x}^{(1)}_k,\tilde{x}^{(2)}_k)$. Let $\tilde{\alpha}_{\mathsf{sca},k}^{(1)}=1-\tilde{\mathbb{P}}^{\gamma,\varepsilon_k}_{N_k}(\tilde{x}_k^{(1)})$ and $\tilde{\alpha}_{\mathsf{sca},k}^{(2)}=1-\tilde{\mathbb{P}}^{\gamma,\varepsilon_k}_{N_k}(\tilde{x}_k^{(2)})$ be sample-based smooth approximation of the violation probability. By repeating the proof of Corollary \ref{coro:z_m_opt}, we can obtain that $\tilde{x}_k^{(1)}\in \tilde{X}_{\tilde{\alpha}_{\mathsf{sca},k}^{(1)}}^{\gamma,\varepsilon_k}(\mathcal{D}_{N_k})$ and $\tilde{x}_k^{(2)}\in \tilde{X}_{\tilde{\alpha}_{\mathsf{sca},k}^{(2)}}^{\gamma,\varepsilon_k}(\mathcal{D}_{N_k})$. Thus, $J(\tilde{x}_k^{(1)})=\tilde{J}^{\gamma,\varepsilon_k}_{\tilde{\alpha}_{\mathsf{sca},k}^{(1)}}(\mathcal{D}_{N_k})$ and $J(\tilde{x}_k^{(2)})=\tilde{J}^{\gamma,\varepsilon_k}_{\tilde{\alpha}_{\mathsf{sca},k}^{(2)}}(\mathcal{D}_{N_k})$. For every $k$, we have
\begin{equation*}
    \tilde{\mathbb{P}}^{\gamma,\varepsilon_k}_{N_k}(\tilde{x}_k^{(1)})\tilde{\mu}_{\mathsf{m}}^{k(1)}+\tilde{\mathbb{P}}^{\gamma,\varepsilon_k}_{N_k}(\tilde{x}_k^{(2)})\tilde{\mu}_{\mathsf{m}}^{k(2)} \geq 1-\alpha.
\end{equation*}
Let $\tilde{\bm{z}}^{\mathsf{lim}}_{\mathsf{m},\alpha}=(\tilde{\mu}_{\mathsf{m}}^{\mathsf{lim}(1)},\tilde{\mu}_{\mathsf{m}}^{\mathsf{lim}(2)},\tilde{x}^{(1)}_{\mathsf{lim}},\tilde{x}^{(2)}_{\mathsf{lim}})$ be an arbitrary cluster point of $\{\tilde{\bm{z}}^k_{\mathsf{m},\alpha}\}_{k=1}^{\infty}$ and $\{\tilde{\bm{z}}^l_{\mathsf{m},\alpha}\}_{l=1}^{\infty}$ be a subsequence that converges to $\tilde{\bm{z}}^{\mathsf{lim}}_{\mathsf{m},\alpha}$. Because $\tilde{\mathbb{P}}^{\gamma,\varepsilon_l}_{N_l}(x)$ converges uniformly to $\mathbb{P}(x)$ on $\mathcal{X}$ w. p. 1 by \cite[Theorem 3.4]{Pena:2020}, we have
\begin{align*}
    \mathbb{P}(\tilde{x}_{\mathsf{lim}}^{(1)})&=\lim_{l\rightarrow\infty}\tilde{\mathbb{P}}^{\gamma,\varepsilon_l}_{N_l}(\tilde{x}_{k}^{(1)})=\tilde{\mathbb{P}}^{\gamma,\varepsilon_{\mathsf{lim}}}_{N_{\mathsf{lim}}}(\tilde{x}_{\mathsf{lim}}^{(1)}), \\
    \mathbb{P}(\tilde{x}_{\mathsf{lim}}^{(2)})&=\lim_{l\rightarrow\infty}\tilde{\mathbb{P}}^{\gamma,\varepsilon_l}_{N_l}(\tilde{x}_{k}^{(2)})=\tilde{\mathbb{P}}^{\gamma,\varepsilon_{\mathsf{lim}}}_{N_{\mathsf{lim}}}(\tilde{x}_{\mathsf{lim}}^{(2)}).
\end{align*}
Then, w. p. 1, we have
\begin{equation*}
    \sum_{i=1}^2\mathbb{P}(\tilde{x}_{\mathsf{lim}}^{(i)})\tilde{\mu}^{\mathsf{lim}(i)}_{\mathsf{m},\alpha}=\sum_{i=1}^2 \tilde{\mathbb{P}}^{\gamma,\varepsilon_{\mathsf{lim}}}_{N_{\mathsf{lim}}}(\tilde{x}_{\mathsf{lim}}^{(i)})\tilde{\mu}^{\mathsf{lim}(i)}_{\mathsf{m},\alpha}\geq 1-\alpha.
\end{equation*}
Therefore, w. p. 1, $\tilde{\bm{z}}^{\mathsf{lim}}_{\mathsf{m},\alpha}$ is a feasible solution of Problem \ref{eq:W_alpha}. Notice that $\tilde{\bm{z}}^l_{\mathsf{m},\alpha}$ converges to $\tilde{\bm{z}}^{\mathsf{lim}}_{\mathsf{m},\alpha}$. Consequently, $\tilde{\mathcal{J}}^{\gamma,\varepsilon_l}_\alpha(\mathcal{D}_{N_l})$ converges to the objective value with $\bm{z}_{\mathsf{m}}=\tilde{\bm{z}}^{\mathsf{lim}}_{\mathsf{m},\alpha}$ which is not smaller than $\mathcal{J}^*_\alpha$. Note that $\tilde{\bm{z}}^{\mathsf{lim}}_{\mathsf{m},\alpha}$ is an arbitrary cluster point. Then, w. p. 1, the following holds:
\begin{equation}
    \label{eq:tilde_calJ_calJ_opt_1}    \liminf_{k\rightarrow\infty}\tilde{\mathcal{J}}^{\gamma,\varepsilon_k}_\alpha(\mathcal{D}_{N_k})=\sum_{i=1}^2 J(\tilde{x}^{(i)}_{\mathsf{lim}})\tilde{\mu}^{\mathsf{lim}(i)}_{\mathsf{m},\alpha}\geq \mathcal{J}^*_\alpha.
\end{equation}
Let $\bm{z}^*_{\mathsf{m},\alpha}=(\mu_{\mathsf{m},\alpha}^{*(1)},\mu_{\mathsf{m},\alpha}^{*(2)},x^{*(1)}, x^{*(2)})\in D_\alpha$ be an optimal solution of Problem \ref{eq:W_alpha}. Let $\tilde{\alpha}_*^{(1)}=1-\mathbb{P}(x^{*(1)})$ and $\tilde{\alpha}_*^{(2)}=1-\mathbb{P}(x^{*(2)})$ be the violation probabilities of $x^{*(1)}$ and $x^{*(2)}$, respectively. By Assumption \ref{assump:opti_solu_Q_alpha}, there exists two sequences $\{x_{l}^{(1)}\}_{l=1}^\infty$ and $\{x_{l}^{(2)}\}_{l=1}^\infty$ converging to $\bm{z}^*_{\mathsf{m}}$ with $\mathbb{P}(x_{l}^{(1)})>1-\tilde{\alpha}_*^{(1)}$ and $\mathbb{P}(x_{l}^{(2)})>1-\tilde{\alpha}_*^{(2)}$, respectively. Then, we have  
\begin{equation}\label{eq:x_l_1_2_P}
    \sum_{i=1}^2 \mathbb{P}(x_{l}^{(i)})\mu_{\mathsf{m},\alpha}^{*(i)}>\sum_{i=1}^2 \mathbb{P}(x_{*}^{(i)})\mu_{\mathsf{m},\alpha}^{*(i)}\geq 1-\alpha.
\end{equation}
Construct $\bm{z}^l_{\mathsf{m},\alpha} :=(\mu_{\mathsf{m},\alpha}^{*(1)},\mu_{\mathsf{m},\alpha}^{*(2)},x^{(1)}_l$ and $x^{(2)}_l$). By \eqref{eq:x_l_1_2_P}, $\bm{z}^l_{\mathsf{m},\alpha}$ is a feasible solution of Problem \ref{eq:W_alpha}. Besides, the sequence $\{\bm{z}^l_{\mathsf{m},\alpha}\}_{l=1}^\infty$ converges to $\bm{z}^*_{\mathsf{m},\alpha}$. Due to \cite[Theorem 3.4]{Pena:2020}, for every $l$ (equivalently pair of $x^{(1)}_{l},x^{(2)}_l$), there exists $K(l)$ such that $\tilde{P}^{\gamma,\varepsilon_k}_{N_k}(x^{(1)}_{l})\geq 1-\tilde{\alpha}^{(1)}_{*}$ and $\tilde{P}^{\gamma,\varepsilon_k}_{N_k}(x^{(2)}_{l})\geq 1-\tilde{\alpha}^{(2)}_{*}$ both hold when $k>K(l)$ w. p. 1. Without loss of generality, assume that $K(l+1)>K(l)$ for every $l$. Then, define a sequence $\{\bar{\bm{z}}^k_{\mathsf{m},\alpha}\}_{k={K(1)}}^{\infty}$ by setting $\bar{\bm{z}}^k_{\mathsf{m},\alpha}=\bm{z}^l_{\mathsf{m},\alpha}$ if $k\in[K(l),K(l+1))$. We use $\bar{\bm{z}}^k_{\mathsf{m},\alpha}=(\mu_{\mathsf{m},\alpha}^{*(1)},\mu_{\mathsf{m},\alpha}^{*(2)},\bar{x}^{(1)}_k,\bar{x}^{(2)}_k)$ with $\bar{x}^{(1)}_k=x^{(1)}_l,\bar{x}^{(2)}_k=x^{(2)}_l$ if $k\in[K(l),K(l+1))$. Then, $\tilde{P}^{\gamma,\varepsilon_k}_{N_k}(\bar{x}^{(1)}_{k})\geq 1-\tilde{\alpha}^{(1)}_{*}$ and $\tilde{P}^{\gamma,\varepsilon_k}_{N_k}(\bar{x}^{(2)}_{k})\geq 1-\tilde{\alpha}^{(2)}_{*}$ holds for every $k$. Thus, the following inequality holds:
\begin{equation}
    \label{eq:fea_tild_x_1_2}
    \sum_{i=1}^2 \tilde{P}^{\gamma,\varepsilon_k}_{N_k}(\bar{x}_{k}^{(i)})\mu_{\mathsf{m},\alpha}^{*(i)} \geq 1- \sum_{i=1}^2 \tilde{\alpha}^{(i)}_{*}\mu_{\mathsf{m},\alpha}^{*(i)}\geq 1-\alpha.
\end{equation}
By \eqref{eq:fea_tild_x_1_2}, we know that $\bar{\bm{z}}_{\mathsf{m},\alpha}$ is a feasible solution of $\tilde{W}_{\alpha}^{\gamma,\varepsilon_k}(\mathcal{D}_{N_k})$ and
\begin{equation*}
    \sum_{i=1}^2 J(\bar{x}^{(i)}_k)\mu_{\mathsf{m},\alpha}^{*(i)} \geq \tilde{\mathcal{J}}_{\alpha}^{\gamma,\varepsilon_k}(\mathcal{D}_{N_k})
\end{equation*}
holds for all $k\geq K(1)$. Since $J(x)$ is continuous due to Assumption \ref{assump:J_h} and $\bar{z}_{\mathsf{m},\alpha}^k$ converges to $z^*_{\mathsf{m},\alpha}$, we have
\begin{equation}\label{eq:tilde_calJ_calJ_opt_2}   
    \limsup_{k\rightarrow\infty}\tilde{\mathcal{J}}^{\gamma,\varepsilon_k}_\alpha(\mathcal{D}_{N_k})\leq \sum_{i=1}^2 J(\bar{x}^{*(i)}) \mu^*_{\mathsf{m},\alpha}(i) =\mathcal{J}^*_\alpha.
\end{equation}
With \eqref{eq:tilde_calJ_calJ_opt_1} and \eqref{eq:tilde_calJ_calJ_opt_2}, we have $\tilde{\mathcal{J}}_{\alpha}^{\gamma,\varepsilon}(\mathcal{D}_N)\rightarrow\mathcal{J}^*_{\alpha}$ w. p. 1. 

We argue by a contradiction for $\mathbb{D}\left(\tilde{D}_{\alpha}^{\gamma,\varepsilon}(\mathcal{D}_N),D_\alpha\right)\rightarrow 0.$ Suppose that $\mathbb{D}\left(\tilde{D}_{\alpha}^{\gamma,\varepsilon}(\mathcal{D}_N),D_\alpha\right)\not\to 0,\ \text{w.p.1.}$ 
Assume there exists $\tilde{\bm{z}}_{\mathsf{m},\alpha}^k\in\tilde{D}_{\alpha}^{\gamma,\varepsilon_k}(\mathcal{D}_{N_k})$ corresponding to $N_k$ and $\varepsilon_k$ in $\{N_k\}_{k=1}^\infty\geq N_0$ and $\{\varepsilon_k\}_{k=1}^{\infty}$ for every $k$, such that $\mathsf{dist}(\tilde{\bm{z}}_{\mathsf{m},\alpha}^k,D_\alpha)\geq \eta$ for some $\eta$. Here, $\mathsf{dist}(\tilde{\bm{z}}_{\mathsf{m},\alpha}^k,D_\alpha)=\min{\|\tilde{\bm{z}}_{\mathsf{m},\alpha}^k-z_{\mathsf{m}}\|:z_{\mathsf{m}}\in D_\alpha}$ is the distance from $\tilde{\bm{z}}_{\mathsf{m},\alpha}^k$ to $D_\alpha$. The limit of $\tilde{\bm{z}}_{\mathsf{m},\alpha}^k$ as $k\rightarrow \infty$ is $\bm{z}^+_{\mathsf{m},\alpha}:=(\mu_{\mathsf{m},\alpha}^{+(1)}, \mu_{\mathsf{m},\alpha}^{+(2)},x^{+(1)},x^{+(2)}),$ which is a feasible solution of \ref{eq:W_alpha} but not in $D_\alpha$. Thus, we have
\begin{equation*}   
    \lim_{k\rightarrow\infty}\tilde{\mathcal{J}}_{\alpha}^{\gamma,\varepsilon}(\mathcal{D}_N) = \sum_{i=1}^2\mu_{\mathsf{m},\alpha}^{+(i)} J(x^{+(i)})>\mathcal{J}^*_\alpha,
\end{equation*}
which contradicts to $\tilde{\mathcal{J}}_{\alpha}^{\gamma,\varepsilon}(\mathcal{D}_N)\rightarrow\mathcal{J}^*_{\alpha}$ w. p. 1. \qed

\section{Proof of Theorem \ref{theo:prob_fea_gua}}
\label{proof:theo_prob_fea_gua}

We briefly summarize Hoeffding's inequality \cite{Hoeffding}, which is applied to prove Theorem \ref{theo:prob_fea_gua}. Let $Y_1,...,Y_N$ be independent random variables, with $\mathsf{Pr}\{Y_j\in[y_{j}^{\mathsf{min}},y_{j}^{\mathsf{max}}]\}=1$, where $y_{j}^{\mathsf{min}}\leq y_{j}^{\mathsf{max}}$ for $j\in [N]$. Then, let $e_i^Y:=Y_i-\mathbb{E}[Y_i]$ and $b_i^Y:=y_{j}^{\mathsf{max}}-y_{j}^{\mathsf{min}}$, if $r>0$, the following condition holds:
\begin{equation}
    \label{eq:Hoeffding_inequa}
    \mathsf{Pr}\left\{\sum_{i=1}^N e_i^Y\geq rN\right\}\leq\exp\left\{-\frac{2N^2r^2}{\sum_{i=1}^N(b_i^Y)^2}\right\}.
\end{equation}

Let $\hat{\bm{z}}_{\mathsf{m}}\in\mathcal{Z}_{\mathsf{m}}$ be a solution out of the feasible set of Problem \ref{eq:W_alpha}, namely, $\hat{\bm{z}}_{\mathsf{m}}\notin\mathcal{Z}_{\mathsf{m},\alpha}$. Specify $\hat{\bm{z}}_{\mathsf{m}}$ by $\hat{\bm{z}}_{\mathsf{m}}=(\hat{\mu}_{\mathsf{m}}^{(1)},\hat{\mu}_{\mathsf{m}}^{(2)},\hat{x}^{(1)},\hat{x}^{(2)})$. 

Let $Y_j=\displaystyle \sum_{i=1}^2\Lambda_{\varepsilon}\left(h(\hat{x}^{(i)},\xi^{(j)})+\gamma\right) \hat{\mu}_{\mathsf{m}}^{(i)}$, $\forall j=1,...,N$, then we have $\mathsf{Pr}\{Y_j\in[0,1]\}=1$ and $\mathbb{E}[Y_j]=\displaystyle \sum_{i=1}^{2}\hat{\mathbb{P}}^{\gamma,\varepsilon}(\hat{x}^{(i)}) \hat{\mu}_{\mathsf{m}}^{(i)}$. Since $\hat{\bm{z}}_{\mathsf{m}}\notin \mathcal{Z}_{\mathsf{m},\alpha}$, we have $\mathbb{P}_{\mathsf{m}}(\hat{\bm{z}}_{\mathsf{m}}):=\sum_{i=1}^{2}\mathbb{P}(\hat{x}^{(i)}) \hat{\mu}_{\mathsf{m}}^{(i)} <1-\alpha.$ 

Define $\tilde{\mathbb{P}}_{\mathsf{m}}(\hat{\bm{z}}_{\mathsf{m}}):=\sum_{i=1}^2 \tilde{\mathbb{P}}^{\gamma,\varepsilon}_N(\hat{x}^{(i)}) \hat{\mu}_{\mathsf{m}}^{(i)},\ \hat{\mathbb{P}}_{\mathsf{m}}(\hat{\bm{z}}_{\mathsf{m}}):=\sum_{i=1}^2 \hat{\mathbb{P}}^{\gamma,\varepsilon}(\hat{x}^{(i)}) \hat{\mu}_{\mathsf{m}}^{(i)}.$
Then, it holds that
\begin{equation}
    \label{eq:Prob_inequality_1_1}
    \tilde{\mathbb{P}}_{\mathsf{m}}(\hat{\bm{z}}_{\mathsf{m}})\geq 1-\alpha'.
\end{equation}
Adding $\hat{\mathbb{P}}_{\mathsf{m}}(\hat{\bm{z}}_{\mathsf{m}})$ to both sides of \eqref{eq:Prob_inequality_1_1}, we obtain
\begin{equation}
    \label{eq:Prob_inequality_1_2}
    \tilde{\mathbb{P}}_{\mathsf{m}}(\hat{\bm{z}}_{\mathsf{m}})-\hat{\mathbb{P}}_{\mathsf{m}}(\hat{\bm{z}}_{\mathsf{m}})\geq 1-\alpha+\alpha-\alpha'-\hat{\mathbb{P}}_{\mathsf{m}}(\hat{\bm{z}}_{\mathsf{m}}).
\end{equation}
Since $\mathbb{P}_{\mathsf{m}}(\hat{\bm{z}}_{\mathsf{m}})<1-\alpha$ holds, using $\mathbb{P}_{\mathsf{m}}(\hat{\bm{z}}_{\mathsf{m}})$ to replace $1-\alpha$ in the left side of \eqref{eq:Prob_inequality_1_2} we have 
\begin{align}
    \label{eq:Prob_inequality_1_3}
    \tilde{\mathbb{P}}_{\mathsf{m}}(\hat{\bm{z}}_{\mathsf{m}})-\hat{\mathbb{P}}_{\mathsf{m}}(\hat{\bm{z}}_{\mathsf{m}})&>\mathbb{P}_{\mathsf{m}}(\hat{\bm{z}}_{\mathsf{m}})-\hat{\mathbb{P}}_{\mathsf{m}}(\hat{\bm{z}}_{\mathsf{m}})+(\alpha-\alpha').
\end{align}

By the definition of $R^{\gamma,\varepsilon}_{\alpha'}$ as \eqref{eq:R_def}, the right side of \eqref{eq:Prob_inequality_1_3} is $R^{\gamma,\varepsilon}_{\alpha'}$, and we have
\begin{align*}
    \tilde{\mathbb{P}}_{\mathsf{m}}(\hat{\bm{z}}_{\mathsf{m}})-\hat{\mathbb{P}}_{\mathsf{m}}(\hat{\bm{z}}_{\mathsf{m}})>R^{\gamma,\varepsilon}_{\alpha'} \label{eq:Prob_inequality_1_4}\ \Rightarrow\ \sum_{j=1}^N (Y_j-\mathbb{E}[Y_j])\geq R^{\gamma,\varepsilon}_{\alpha'}N.  
\end{align*}
Therefore, we have
\begin{align}
    \mathsf{Pr}\{\hat{\bm{z}}_{\mathsf{m}}\in\mathcal{Z}^{\gamma,\varepsilon}_{\mathsf{m},\alpha'}(\mathcal{D}_N)\}&=\mathsf{Pr}\{\tilde{\mathbb{P}}_{\mathsf{m}}(\hat{\bm{z}}_{\mathsf{m}})\geq 1-\alpha'\} \nonumber \\
    &\leq\mathsf{Pr}\{\sum_{j=1}^N (Y_j-\mathbb{E}[Y_j])\geq R^{\gamma,\varepsilon}_{\alpha'}N\} \nonumber\\
    &\leq\exp\{-2N\left(R^{\gamma,\varepsilon}_{\alpha'}\right)^2\}, \nonumber
\end{align}  
where the last inequality holds due to Hoeffding's inequality \eqref{eq:Hoeffding_inequa}. Recall that $\tilde{\bm{z}}_{\mathsf{m},\alpha'}$ is an optimal solution of Problem \ref{eq:til_W_alpha_gamma_N}. If we suppose that $\tilde{\bm{z}}_{\mathsf{m},\alpha'}$ is not feasible for Problem \ref{eq:W_alpha}, then the following holds
\begin{equation}
\label{eq:z_m_alpha_ineq_1}
    \tilde{\bm{z}}_{\mathsf{m},\alpha'}\in \tilde{D}^{\gamma,\varepsilon}_{\alpha'}(\mathcal{D}_N),\ \tilde{\bm{z}}_{\mathsf{m},\alpha'}\notin \mathcal{Z}_{\mathsf{m},\alpha}.
\end{equation}
Notice that \eqref{eq:z_m_alpha_ineq_1} implies that
\begin{equation}
    \label{eq:z_m_alpha_ineq_2}
    \tilde{\bm{z}}_{\mathsf{m},\alpha'}\in\mathcal{Z}^{\gamma,\varepsilon}_{\mathsf{m},\alpha'}(\mathcal{D}_N),\ \tilde{\bm{z}}_{\mathsf{m},\alpha'}\notin \mathcal{Z}_{\mathsf{m},\alpha}.
\end{equation}

From \eqref{eq:z_m_alpha_ineq_2}, we further have
\begin{equation}
    \label{eq:z_m_alpha_ineq_3}
    \hat{\bm{z}}_{\mathsf{m}}\in\mathcal{Z}^{\gamma,\varepsilon}_{\mathsf{m},\alpha'}(\mathcal{D}_N), \hat{\bm{z}}_{\mathsf{m}}\notin \mathcal{Z}_{\mathsf{m},\alpha}.
\end{equation}
Finally, we conclude 
\begin{align*}
    \mathsf{Pr}\{\tilde{\bm{z}}_{\mathsf{m},\alpha'}\notin\mathcal{Z}_{\mathsf{m},\alpha}\}&\leq\mathsf{Pr}\{\tilde{\bm{z}}_{\mathsf{m},\alpha'}\in\mathcal{Z}^{\gamma,\varepsilon}_{\mathsf{m},\alpha'}(\mathcal{D}_N)\}
    \nonumber \\
    &\leq\exp\{-2N\left(R^{\gamma,\varepsilon}_{\alpha'}\right)^2\}. \qed
\end{align*} 

\end{document}